\documentclass[12pt,reqno]{amsart}
\usepackage{fullpage}
\usepackage{hyperref}
\usepackage{amsrefs}

\newtheorem*{thm*}{Theorem}
\newtheorem{thm}{Theorem}[section]
\newtheorem{prop}[thm]{Proposition}
    \newtheorem{dfn}[thm]{Definition}
\theoremstyle{remark}
    
    \newtheorem*{rem*}{Remark}
\numberwithin{equation}{section}

\newcommand{\mcJ}{\mathcal{J}}
\newcommand{\fg}{\mathfrak g}
\newcommand{\wF}[2]{\widetilde F_{(#2)}^{#1}}
\newcommand{\wE}[2]{\widetilde E_{(#2)}^{#1}}
\newcommand{\wV}[2]{\widetilde V_{(#2)}^{#1}}
\renewcommand{\d}{\partial}
\newcommand{\TT}{T\!\!\!T}
\newcommand{\gl}{{\mathfrak g}{\mathfrak l}}
\newcommand{\so}{{\mathfrak s}{\mathfrak o}}
\newcommand{\spm}{{\mathfrak s}{\mathfrak p}}
\newcommand{\Har}{{\mathcal H}}
\newcommand{\Gr}{{\mbox{Gr}_{\widehat A^+}}}
\newcommand{\Hom}{{\mbox{Hom}}}
\newcommand{\rank}{{\mbox{rank}}}
\newcommand{\nt}{{\noindent}}
\newcommand{\C}{{{\mathbb C}}}

\newcommand{\calS}{{\mathcal S}}
\newcommand{\calp}{{\mathcal P}}
\newcommand{\calr}{{\mathcal R}}
\newcommand{\calh}{{\mathcal H}}

\newcommand{\bbC}{\mathbb{C}}
\newcommand{\mcP}{\mathcal{P}}
\newcommand{\mcD}{\mathcal{D}}
\newcommand{\mcH}{\mathcal{H}}
\newcommand{\mcS}{\mathcal{S}}
\newcommand{\mcW}{\mathcal{W}}
\newcommand{\la}{\lambda}
\newcommand{\al}{\alpha}
\renewcommand{\sp}{\mathfrak{sp}}
\newcommand{\F}[2]{F^{#1}_{(#2)}}
\newcommand{\E}[2]{E^{#1}_{(#2)}}
\newcommand{\V}[2]{V^{#1}_{(#2)}}
\newcommand{\ot}{\otimes}
\newcommand{\dual}[1]{ {\left({#1}\right)}^*}

\begin{document}

\title{Reciprocity Algebras and Branching for Classical Symmetric Pairs}
\date{\today}

\author{Roger E. Howe}
\author{Eng-Chye Tan}
\author{Jeb F. Willenbring}
\address{Department of Mathematics, Yale University, Box 208283 Yale Station, New Haven, CT 06520, USA}
\address{Department of Mathematics, National University of Singapore, 2 Science Drive 2, Singapore 117543, Singapore}
\address{Department of Mathematical Sciences, University of Wisconsin at Milwaukee, Milwaukee, WI 53201, USA}

\begin{abstract}
We study branching laws for a classical group $G$ and a symmetric
subgroup $H$.  Our approach is through the {\it branching algebra},
the algebra of covariants for $H$ in the regular functions on the
natural torus bundle over the flag manifold for $G$.  We give
concrete descriptions of certain subalgebras of the branching
algebra using classical invariant theory.  In this context, it turns
out that the ten classes of classical symmetric pairs $(G,H)$ are
associated in pairs, $(G,H)$ and $(H',G')$, and that the (partial)
branching algebra for $(G,H)$ also describes a branching law from
$H'$ to $G'$.  (However, the second branching law may involve
certain infinite-dimensional highest weight modules for $H'$.)  To
highlight the fact that these algebras describe two branching laws
simultaneously, we call them {\it reciprocity algebras}.  Our
description of the reciprocity algebras reveals that they all are
related to the tensor product algebra for $GL_n$.  This relation is
especially strong in the {\it stable range}.  We provide explicit
descriptions of reciprocity algebras in the stable range in terms of
the tensor product algebra for $GL_n$.  This is the structure lying
behind formulas for branching multiplicities in terms of
Littlewood-Richardson coefficients.
\end{abstract}

\maketitle

\section{Introduction}

One of the basic problems in the study of transformation groups
and finite dimensional representation theory is the \it branching
problem\rm: describing how irreducible representations of a
compact group (\underbar {mutatis} \underbar {mutandis}, a
reductive complex algebraic group) decompose when restricted to a
subgroup. More specifically, given a group $G$ and a subgroup $H$,
we may consider the branching problem for the pair $(G, H)$. This
paper introduces and establishes some basic features of an
approach to the branching problem for a class of particularly
interesting cases: the classical symmetric pairs, in which the
group $G$ is a classical group and the subgroup $H$ is a symmetric
subgroup - the fixed points of an involution of $G$. One can
organize these pairs in 10 families, listed in table I (see page
19).

There is already a substantial literature on branching rules for
classical symmetric pairs \cite{benkart},
\cite{berenstein-zelevinsky}, \cite{berenstein-zelevinsky-2},
\cite{black-king-wybourne}, \cite{ew-annals}, \cite{fulton},
\cite{goodman-wallach-book}, \cite{howe-reciprocity},
\cite{howe-schur}, \cite{HTW1}, \cite{HTW2}, \cite{james-kerber},
\cite{ king-modify}, \cite{king-tensor}, \cite{king-plethysm},
\cite{king-spinor}, \cite{king-s-fcn}, \cite{king-affine}, \cite{
king-analogies-I}, \cite{king-analogies-II}, \cite{knapp-branching},
\cite{knapp-book}, \cite{ koike-terada-littlewood},
\cite{koike-terada-Sp-SO}, \cite{ koike-terada-classical},
\cite{koike-terada}, \cite{kostant-new}, \cite{kudla-see-saw},
\cite{littelmann-1}, \cite{littelmann-2}, \cite{littelmann-3},
\cite{littelmann}, \cite{littlewood-richardson}, \cite{littlewood},
\cite{littlewood-paper}, \cite{newell}, \cite{popov}, \cite{rees},
\cite{sundaram}, \cite{vinberg}, \\ \cite{weyl}, \cite{zhelobenko}.

A substantial portion of this work is combinatorial in nature, and
is formulated in terms of tableaux and Littlewood-Richardson
coefficients. This approach began with the work of Littlewood and
Richardson \cite{littlewood-richardson} \cite{fulton} on tensor
products of representations of $GL_n$ and Littlewood's Restriction
Formula \cite{littlewood} \cite{littlewood-paper} for restriction
of certain representations of $GL_n$ to the orthogonal group or
the symplectic group. More recent work in this tradition includes
\cite{sundaram}, several papers by Koike and Terada
\cite{koike-terada-littlewood}, \cite{koike-terada-Sp-SO},
\cite{koike-terada-classical}, \cite{koike-terada} and by R.W.
King and collaborators \cite{king-modify}, \cite{king-tensor},
\cite{king-plethysm}, \cite{king-spinor},
\cite{black-king-wybourne}, \cite{king-s-fcn},
\cite{king-analogies-I}, \cite{king-analogies-II},
\cite{king-affine}.

A group $G$ embedded diagonally in the product $G \times G$ is a
symmetric subgroup, corresponding to the involution which exchanges
the two copies of $G$. An irreducible representation of $G \times G$
is just a tensor product $\rho_1 \otimes \rho_2$ of irreducible
representations of $G$. Thus, the branching problem for symmetric
pairs includes the \it tensor product problem \rm - to describe how
a tensor product of representations of a group $G$ decomposes into
irreducible pieces. In the 1990s, striking progress on the tensor
product problem came from an unexpected source -- the theory of
quantum groups (see \cite{drinfeld, jantzen, jimbo, joseph}).

Quantum groups are one-parameter deformations of the universal
enveloping algebra of a semisimple Lie algebra. The extra parameter
was used by Lusztig (\cite{lusztig}) and Kashiwara
(\cite{kashiwara}) to define bases (\lq \lq canonical bases" or \lq
\lq crystal bases") for these algebras. The properties of these
bases permitted combinatorial descriptions of tensor product
multiplicities, providing an extension of the Littlewood-Richardson
Rule for $GL_n$ to all simple groups. This approach has been
elucidated and simplified by the path model of Littelmann (see
\cite{littelmann-1, littelmann-2, littelmann-3, littelmann}), which
also provides branching rules for restriction from a group $G$ to
the Levi component of a parabolic subgroup.

This paper takes another approach to the branching problem for
classical symmetric pairs. Its main tool is invariant theory,
represented primarily by the geometric version of highest weight
theory (see \cite{goodman-wallach-book}, \cite{howe-schur} and
\cite{bott}) and by classical invariant theory \cite{weyl,
howe-schur, goodman-wallach-book}.  (At the moment, parts of our
development also depend on the combinatorial theory of the
Littlewood-Richardson Rule. However, our approach also suggests
alternative approaches to these results, and we hope eventually to
obtain independent proofs.) To describe the branching rule for a
classical symmetric pair $(G, H)$, we study a naturally defined
algebra $B(G, H)$, which we call the \it branching algebra \rm for
the pair $(G, H)$. The algebra $B(G, H)$ carries a multigrading by
$\hat A_G^+ \times \hat A_H^+$, the semigroups of dominant weights
of $G$ and $H$. The dimension of the homogeneous component $B(G,
H)^{(\phi,\psi)}$ labeled by $(\phi, \psi) \in \hat A_G^+ \times
\hat A_H^+$ equals the multiplicity of $\rho^{\psi}$ in
$\rho^{\phi}$, where $\rho^{\phi}$ (respectively, $\rho^{\psi}$)
is the irreducible representation of $G$ with highest weight
$\phi$ (respectively, irreducible representation of $H$, with
highest weight $\psi$). Branching algebras were studied in
\cite{zhelobenko}, but have been mostly neglected since.

In addition to the full branching algebras, our theory distinguishes
certain families of graded subalgebras, depending on integer
parameters. For a certain range of the parameters, the subalgebras
have a particularly simple structure. We call this the \it stable
range, \rm and the resulting algebras, \it stable branching
algebras. \rm

\medskip
This approach must of course give the same numbers as are produced
by the combinatorial, quantum group and path model methods. We have
already shown in \cite{HTW1} how to use branching algebras
systematically to express stable branching multiplicities for all
families of classical symmetric pairs in terms of
Littlewood-Richardson coefficients.

We feel that branching algebras have several virtues in addition
to providing an efficient route to formulas for branching
multiplicities. One is, they provide finer information than
multiplicities: they provide highest weight vectors. Indeed, the
homogeneous component $B(G, H)^{(\phi, \psi)}$ of $B(G,H)$
consists of the $\psi$-highest weight vectors for $H$ in a
specific realization of the representation $\rho^{\phi}$ of $G$.
Thus, if one understands $B(G,H)$ well enough, one gets
information about how $H$-modules sit inside $G$-modules. For the
classical symmetric pairs, this possibility can in fact be
realized: in \cite{HTW2}, the authors have shown how to construct
an explicit basis for the $GL_n$ tensor product algebra. In
further work, Howe and Lee \cite{howe-lee} have extended this
construction to apply to most of the families of branching
algebras.

A second positive feature of branching algebras is, they collect all
multiplicity information into one coherent structure. This has
several potential implications. For example, it has been observed in
many cases that branching multiplicities can be interpreted as the
number of integral points in certain convex polytopes. The theory of
SAGBI bases (\cite{SAGBI}) provides a strategy for attaching
lattices cones to affine algebras. The explicit bases mentioned
above permit one to realize the SAGBI strategy effectively. The
insight gained into the structure of branching algebras allows one
in turn to predict the existence of convex polytopes whose integral
points define the desired multiplicities
\cite{berenstein-zelevinsky}, \cite{berenstein-zelevinsky-2},
\cite{lak-1}, \cite{lak-2}.

The algebra structure also provides natural connections between
the different branching algebras. The stable branching algebras
for all 10 families of classical symmetric pairs are closely
related to each other, and in particular, they all can be
described in terms of the $GL_n$ branching algebra. These
relationships lie behind the various formulas expressing stable
multiplicities for all families in terms of Littlewood-Richardson
coefficients \cite{HTW1}.

These connections in fact go beyond branching algebras. By some
remarkable coincidences, there are also very close relationships
between branching algebras for classical symmetric pairs and
classical Kostant-Rallis actions (see \cite{kostant-rallis},
\cite{wallach-willenbring}, \cite{tams-willenbring}). By a
Kostant-Rallis action, we mean the complexified action of $H$ on
$\fg / {\mathfrak h}$, where $H$ is a maximal compact subgroup of
a semisimple Lie group $G$, and $\fg$ and $\mathfrak h$ are the
Lie algebras of $G$ and $H$ respectively. Some key results about
these actions were given in \cite{kostant-rallis}. but some
aspects of their structure are still not completely settled. The
tight connection of classical Kostant-Rallis actions with
branching algebras can be used to improve understanding of these
important examples. Beyond these immediate applications. branching
algebras should provide a useful general tool for analyzing group
actions.

We also feel that the branching algebra approach sheds light on some
of the basic constructions of the combinatorial theory. For example,
the combinatorial description tensor product multiplicities for
$GL_n$ by means of Littlewood-Richardson tableaux still has an air
of mystery about it. A number of experts have remarked on the
curious circumstance that, although tensor product multiplicities
are obviously symmetric in the two tensor factors, the
Littlewood-Richardson description of multiplicities is not
symmetric. The SAGBI basis description of the $GL_n$ tensor product
algebra \cite{HTW2} provides a natural interpretation of the
information encoded by a Littlewood-Richardson tableau, including
providing a straightforward context for the asymmetrical treatment
of the two factors.

The efficient packaging of all multiplicity formulas for a given
pair in one algebraic structure also seems to provide an efficient
way of computing explicit results in low rank examples. The
authors hope to treat a collection of explicit cases in a future
publication.

Finally, although the branching algebra construction is not
limited to symmetric pairs, it has a remarkable additional feature
in the case of classical symmetric pairs. As explained in
\cite{howe-reciprocity} \cite{howe-remarks}, classical invariant
theory gives rise to natural correspondences between
representations of pairs (\lq \lq dual pairs") of classical
groups. It turns out that these correspondences can be refined to
correspondences between branching laws for symmetric pairs. The 10
families of classical symmetric pairs pair up into 5 families of
pairs of pairs (see Table II on page 19), in such a way that a
branching algebra for one pair can also be interpreted as a
branching algebra for the corresponding pair. Thus, one branching
algebra simultaneously describes branching between representations
of two classical symmetric pairs. This can be thought of as a
reciprocity law, analogous to Frobenius Reciprocity in the theory
of group representations. Hence we call these branching algebras
\it reciprocity algebras. \rm (Multiplicity versions of these
reciprocity laws are implicitly included in the general result of
\cite{howe-reciprocity}). We should also remark that in most
cases, the relevant representations for some of the groups
involved in the pairs will be infinite dimensional. However, the
infinite dimensional representations which occur are of a very
special sort, and share many of the characteristics of
finite-dimensional representations.)

For all the above reasons, we hope that branching algebras will
serve as a useful complement to the already developed approaches to
describing branching laws.

\medskip
\subsection{Overview}

We will give the construction of a branching algebra in the
following section. We shall illustrate how the algebra structure
emphasizes the coherence of the branching problem across
representations of $G$, rather than simply providing a
representation-by-representation answer. Further, the algebra
structure goes beyond the numerical nature of multiplicities, and
is the basis of some of the well-known multiplicity results.

There are three main parts to this paper:
\begin{enumerate}
\item[(a)] {\bf Background:} preliminary knowledge on
multiplicity-free spaces, the theory of dual pairs and the
classification of classical symmetric pairs, leading to the
concept of reciprocity algebra. One can group the classical
symmetric pairs into 10 families, and it turns out that these
families are associated in pairs, $(G,H)$ and $(H',G')$ (see
Tables I and II), and that the branching algebra for $(G,H)$ also
describes a branching law from $H'$ to $G'$.  Hence, we call the
algebra that describes the two related branching laws a {\it
reciprocity algebra}.

\item[(b)] {\bf Reciprocity:} explaining the underpinnings of the
reciprocity of multiplicities for the 5 pairs of reciprocity
algebras, with detailed discussions for certain representative
pairs.

\item[(c)] {\bf Stability:} describing the branching algebra for a
symmetric pair $(G,H)$ in terms of a $GL_n$ tensor product
algebra. This emphasizes the underlying importance of the $GL_n$
branching algebras, since their structure relate to all others.
This phenomenon is also the basis for the expression of
multiplicity formulas for representations of classical groups in
terms of Littlewood-Richardson coefficients - a recurring theme in
the literature.  The present paper brings more structure to these
multiplicity results.

\end{enumerate}

The first four sections provide the necessary background:
\begin{enumerate}
\item[\S 3:] We provide the notations for the parametrization of
representations in \S 3.1. Next, we discuss some preliminary
concepts on multiplicity-free spaces and prove a main theorem (see
Theorem 3.3) on embedding the associated graded algebra of a
graded $G$-multiplicity-free space into the algebra ${\mathcal R}
(G/U_G)$. The final subsection discusses the rudiments of
classical invariant theory, formulated in terms of the theory of
dual pairs.

\item[\S 4:] Section 4 begins with the simplest example of a pair
of branching algebras for $GL_n$, leading us to the concept of
reciprocity algebras. This is where we bring forth the
classification of classical symmetric pairs and the theory of dual
pairs to introduce the 5 pairs (also called {\it see-saw pairs})
of reciprocity algebras.
\end{enumerate}

For the second part of the  paper, we organize discussions of the
see-saw pairs as follows:
\begin{enumerate}
\item[\S 5:] discusses the branching from $GL_n$ to $O_n$ (note
that branching from $GL_{2n}$ to $Sp_{2n}$ can be treated
similarly);

\item[\S 6:] discusses the tensor product algebra for $O_n$ (note
that the tensor product algebra of $Sp_{2n}$ can be treated
similarly);

\item[\S 7:] provides a more general treatment of the reciprocity
for $(GL_n,GL_m)$.
\end{enumerate}
\noindent For the sake of brevity, we will not discuss the
branching from $Sp_{2(m+n)}$ to $Sp_{2m} \times Sp_{2n}$, $O_{2n}$
to $GL_{n}$ and $Sp_{2n}$ to $GL_n$.  They can be treated
similarly.

In the final part of this paper we begin with some important
comments on stability results for branching in \S 8. We demonstrate
the theoretical underpinnings for stability results, highlighting
certain specific see-saw pairs:

\begin{enumerate}
\item[\S 9:] we interpret the associated graded of the branching
algebra from $GL_n$ to $O_n$ as a $(0,1)$-subalgebra (see
Definition 3.1) of the tensor product algebra of $GL_m$ in the
stable range $n > 2m$.

\item[\S 10:] we interpret the associated graded of the $O_n$
tensor product algebra as a $(0,1)$-subalgebra of a triple product
of tensor product algebras of $GL_n$, $GL_m$ and $GL_{\ell}$ in
the stable range $n> 2(m+\ell)$.

\item[\S 11:] we interpret the branching algebra of $O_{n+m}$ to
$O_n \times O_m$ as a $(0,1)$-subalgebra of a triple tensor
product algebra of $GL_{\ell}$ in the stable range $\min\,(n, m) >
2\ell $.
\end{enumerate}

In all these cases, we show that the branching algebras associated
to symmetric pairs can all be identified with suitable branching
algebras associated to the general linear groups.  Thus, if we can
have control of the solution in the general linear group case, we
will have some control of the other classical groups.  And indeed,
we do (see \cite{HTW2}). The other non-trivial examples will be
important extensions of this work, and we hope to see them in
further papers, for example, \cite{howe-lee}.

\medskip

\nt {\bf Acknowledgements:} We thank Kenji Ueno and Chen-Bo Zhu for
discussions on the proof in the Appendix. The second named author
also acknowledges the support of NUS grant R-146-000-050-112.  The
third named author was supported by NSA Grant \# H98230-05-1-0078.

\section{Branching Algebras}

For a reductive complex linear algebraic $G$, let $U_G$ be a maximal
unipotent subgroup of $G$. The group $U_G$ is determined up to
conjugacy in $G$ \cite{borel-book}. Let $A_G$ denote a maximal torus
which normalizes $U_G$, so that $B_G = A_G \cdot U_G$ is a Borel
subgroup of $G$. Also let $\widehat A_G^+$ be the set of dominant
characters of $A_G$ -- the semigroup of highest weights of
representations of $G$. It is well-known \cite{borel-book}
\cite{howe-schur} and may be thought of as a geometric version of
the theory of the highest weight, that the space of regular
functions on the coset space  $G/U_G$, denoted by $\calr (G/U_G)$,
decomposes (under the action of $G$ by left translations) as a
direct sum of one copy of each irreducible representation $V_{\psi}$
(with highest weight $\psi$) of $G$ (see \cite{towber}):
$$
\calr (G/U_G) \simeq \bigoplus_{ \psi \in \widehat A_G^+}
V_{\psi}. \eqno (2.1)
$$

We note that $\calr (G/U_G)$
has the structure of an $\widehat A_G^+$-graded algebra, for which the $V_{\psi}$ are
the graded components. To be specific, we note that since $A_G$ normalizes
$U_G$, it acts on $G/U_G$ by right translations, and this action commutes
with the action of $G$ by left translations.  The following result is well-known, probably folklore. We provide the proof anyway, since we are not able to find a suitable reference.

\begin{prop}  The algebra of regular functions $\calr (G/U_G)$ is an
$\widehat A_G^+$-graded algebra, under the right action of $A_G$.
More precisely, the decomposition (2.1) is the graded algebra
decomposition under $A_G$, where $V_\psi$ is the $A_G$-eigenspace
corresponding to $\phi \in \widehat A_G^+$ with $\phi =
w^*(\psi^{-1})$.  Here $w$ is the longest element of the Weyl
group with respect to the root system determined by the Borel
subgroup $B_G$.
\end{prop}

\nt {\bf Proof.} Since $A_G$ is commutative and reductive, we can
decompose $\calr (G/U_G)$ into eigenspaces for the right action of
$A_G$. These eigenspaces must be invariant under the action of $G$
by left translations. Since $V_{\psi}$ is irreducible for the
action of $G$, it must belong to a single $A_G$ eigenspace. Let
$f$ be the highest weight vector for $B_G$ in $V_{\psi}$. Then by
definition $L_a(f) = \psi(a) (f)$, where $a \in A_G$ and $L_g$
indicates left translation by $g$: $L_g(f) (h) = f(g^{-1}h)$ for
any $h$ in $G$, representing a point in $G/U_G$. If $w$ is the
longest element of the Weyl group, then the $B_G$ orbit $B_GwU_G$
is dense in $G/U_G$, so that $f$ is determined by its restriction
to $B_GwU_G$. We compute that
$$
\psi (a) f (w) = L_a (f) (w) = f(a^{-1}w) = f(w w^{-1} a^{-1}w) =
R_{w^{-1}a^{-1}w}f(w),
$$
where $R_g f(h) = f(hg)$ refers to the right translation of $f$ by
$g$.  If $V_{\psi}$ belongs to the $\phi$ eigenspace for the right
action of $A_G$, then this equation implies that
$$
\psi (a) = \phi(w^{-1} a^{-1} w), \quad \quad \quad {\rm or }
\quad \quad \quad \psi = w^*(\phi^{-1}), \eqno (2.2)
$$
where $w^*$ indicates the action of $w$ on $\widehat A_G$
resulting from conjugation of $A_G$ by $w$. Thus, the $V_{\psi}$
are exactly the eigenspaces for the right action of $A_G$, with
the $A_G$-eigencharacter related to the highest weight by the
equation (2.2). Since the right action by $A_G$ (as well as the
left action by $G$) is an action by algebra automorphisms of
$\calr (G/U_G)$, it is easy to check that if $f_1$ is a
$\phi$-eigenfunction for $A_G$, and $f_2$ is a
$\theta$-eigenfunction, then the product $f_1f_2$ is a $\phi
\theta$-eigenfunction for $A_G$. It follows that
$$
V_{\psi} V_{\chi} = V_{\psi \chi},
$$
so that, indeed, the decomposition (2.1) defines a structure of an
$\widehat A_G^+$-graded algebra on $\calr (G/U_G).\qquad \qquad
\square$
\medskip

Now let  $H \subset G$ be a reductive subgroup, and let $U_H$ be a
maximal unipotent subgroup of $H$. We consider the algebra $\calr
(G/U_G)^{U_H}$, of functions on $G/U_G$ which are invariant under
left translations by $U_H$. Let $A_H$ be a maximal torus of $H$
normalizing $U_H$, so that $B_H := A_H \cdot U_H$ is a Borel
subgroup of $H$. Then $\calr (G/U_G)^{U_H}$ will be invariant
under the (left) action of $A_H$, and we may decompose $\calr
(G/U_G)^{U_H}$ into eigenspaces for $A_H$. Since the functions in
$\calr (G/U_G)^{U_H}$ are by definition (left) invariant under
$U_H$, the (left) $A_H$-eigenfunctions will in fact be (left)
$B_H$ eigenfunctions. In other words, they are highest weight
vectors for $H$. Hence, the characters of $A_H$ acting on (the
left of) $\calr (G/U_G)^{U_H}$ will all be dominant with respect
to $B_H$, and we may write $\calr (G/U_G)^{U_H}$ as a sum of
(left) $A_H$ eigenspaces $(\calr (G/U_G)^{U_H})^{\chi}$ for
dominant characters $\chi$ of $H$:
$$
\calr (G/U_G)^{U_H} = \bigoplus_{\chi \in \widehat A_H ^+}( \calr
(G/U_G)^{U_H})^{\chi}. \eqno (2.3)
$$

Since the spaces $V_{\psi}$ of decomposition (2.1) are (left)
$G$-invariant, they are \underbar{a} \underbar{fortiori} \it left
\rm $H$-invariant, so we have a decomposition of $\calr
(G/U_G)^{U_H}$ into \it right \rm $A_G$-eigenspaces $( \calr
(G/U_G)^{U_H})_{\psi}$:
$$
\calr (G/U_G)^{U_H} = \bigoplus_{\psi \in \widehat A_G^+} \calr (G/U_G)^{U_H} \cap V_{\psi} :=
\bigoplus_{\psi \in \widehat A_G^+} \calr (G/U_G)^{U_H}_{\psi}.
$$
Combining this decomposition with the decomposition (2.3), we may
write
$$
\calr (G/U_G)^{U_H} = \bigoplus_{\psi \in \widehat A_G^+,\ \chi
\in \widehat A_H ^+} ( \calr (G/U_G)^{U_H}_{\psi})^{\chi}. \eqno
(2.4)
$$

To emphasize the key features of this algebra, we note the
resulting consequences of decomposition (2.4) in the following
proposition.
\begin{prop}
\begin{enumerate}
\item[(a)] The decomposition (2.4) is an ($\widehat A_G^+ \times
\widehat A_H ^+$)-graded algebra decomposition of $\calr
(G/U_G)^{U_H}$. \item[(b)] The subspaces $(\calr
(G/U_G)^{U_H}_{\psi})^{\chi}$ tell us the $\chi$ highest weight
vectors for $B_H$ in the irreducible representation $V_{\psi}$ of
$G$. Therefore, the decomposition
$$
\calr (G/U_G)^{U_H}_{\psi} = \bigoplus_{\chi \in \widehat A_H ^+} ( \calr (G/U_G)^{U_H}_{\psi})^{\chi}
$$
tells us how $V_{\psi}$ decomposes as a $H$-module.
\end{enumerate}
\end{prop}

Thus, knowledge of $\calr (G/U_G)^{U_H}$ as a ($\widehat A_G^+
\times \widehat A_H ^+$)-graded algebra tell us how
representations of $G$ decompose when restricted to $H$, in other
words, it describes the branching rule from $G$ to $H$. We will
call $\calr (G/U_G)^{U_H}$ the $(G, H)$ \it branching algebra. \rm
When $G \simeq H \times H$, and $H$ is embedded diagonally in $G$,
the branching algebra describes the decomposition of tensor
products of representations of $H$, and we then call it the \it
tensor product algebra \rm for $H$.   More generally, we would
like to understand the $(G, H)$ branching algebras for symmetric
pairs $(G, H)$.

\section{Preliminaries and Notations}

\subsection{Parametrization of Representations} Let $G$ be a classical reductive algebraic group
over $\C$: $G = GL_n({\C})=GL_n$, the general linear group; or $G =
O_n ({\C})=O_n$, the orthogonal group; or $G =
Sp_{2n}({\C})=Sp_{2n}$, the symplectic group. We shall explain our
notations on irreducible representations of $G$ using integer
partitions.  In each of these cases, we select a Borel subalgebra of
the classical Lie algebra and coordinatize it, as is done in
\cite{goodman-wallach-book}. Consequently, all highest weights are
parameterized in the standard way (see \cite{goodman-wallach-book}).

A non-negative integer {\it partition} $\lambda$, with $k$ parts, is
an integer sequence $\lambda_1 \geq \lambda_2 \geq \ldots \geq
\lambda_k > 0$.  We may sometimes refer to $\lambda$ as a {\it
Young} or {\it Ferrers diagram}. We use the same notation for
partitions as is done in \cite{macdonald}.  For example, we write
$\ell(\la)$ to denote the {\it length} (or {\it depth}) of a
partition, i.e., $\ell(\la)=k$ for the above partition. Also let
$|\la| = \sum_i \la_i$ be the size of a partition and $\la^\prime$
denote the {\it transpose} (or {\it conjugate}) of $\la$ (i.e.,
$(\la^\prime)_i = |\{\la_j: \la_j \geq i\}|$).

\vskip 10pt

\noindent{\bf ${\bf GL_n}$ Representations:} Given non-negative integers
$p$, $q$ and $n$ such that $n \geq p+q$ and non-negative integer
partitions $\la^+$ and $\la^-$ with $p$ and $q$ parts
respectively, let $\F{(\la^+, \la^-)}{n}$ denote the irreducible
rational representation of $GL_n$ with highest weight given by the
$n$-tuple:
\[
    \begin{array}{ccc}
        (\la^+,\la^-) &=& \underbrace{\left(\la^+_1, \la^+_2, \cdots, \la^+_p, 0, \cdots, 0, -\la^-_{q},
\cdots, -\la^-_{1} \right)} \\
        & & n
    \end{array}
\]
If $\la^- = (0)$ then we will write $\F{\la^+}{n}$ for $\F{(\la^+,
\la^- )}{n}$. Note that if $\la^+ = (0)$ then
$\dual{\F{\la^-}{n}}$ is equivalent to $\F{(\la^+, \la^-)}{n}$.  More generally,
$\dual{\F{(\la^+,\la^-)}{n}}$ is equivalent to $\F{(\la^-, \la^+)}{n}$.

\vskip 10pt

\noindent{\bf ${\bf O_n}$ Representations:} The complex orthogonal
group has two connected components. Because the group is
disconnected we cannot index irreducible representations by highest
weights. There is however an analog of Schur-Weyl duality for the
case of $O_n$ in which each irreducible rational representation is
indexed uniquely by a non-negative integer partition $\nu$ such that
$(\nu^\prime)_1 + (\nu^\prime)_2 \leq n$.  That is, the sum of the
first two columns of the Young diagram of $\nu$ is at most $n$. We
will call such a diagram $O_n$-admissible (see
\cite{goodman-wallach-book} Chapter 10 for details).  Let
$\E{\nu}{n}$ denote the irreducible representation of $O_n$ indexed
$\nu$ in this way.

An irreducible rational representation of $SO_n$ may be indexed by
its highest weight.  In \cite{goodman-wallach-book} Section 5.2.2,
the irreducible representations of $O_n$ are determined in terms of
their restrictions to $SO_n$ (which is a normal subgroup having
index 2). We note that if $\ell (\nu) \neq \frac{n}{2}$, then the
restriction of $\E{\nu}{n}$ to $SO_n$ is irreducible.  If $\ell
(\nu) = \frac{n}{2}$ ($n$ even), then $\E{\nu}{n}$ decomposes into
exactly two irreducible representations of $SO_n$. See
\cite{goodman-wallach-book} Section 10.2.4 and 10.2.5 for the
correspondence between this parametrization and the above
parametrization by partitions.

The determinant defines an (irreducible) one-dimensional
representation of $O_n$.  This representation is indexed by the
length $n$ partition $\zeta = (1,1, \cdots, 1)$.  An irreducible
representation of $O_n$ will remain irreducible when tensored by
$\E{\zeta}{n}$, but the resulting representation \emph{may} be
inequivalent to the initial representation.  We say that a pair of
$O_n$-admissible partitions $\alpha$ and $\beta$ are
\emph{associate} if $\E{\alpha}{n} \otimes \E{\zeta}{n} \cong
\E{\beta}{n}$.  It turns out that $\alpha$ and $\beta$ are associate
exactly when $(\alpha^\prime)_1 + (\beta^\prime)_1 = n$ and
$(\alpha^\prime)_i = (\beta^\prime)_i$ for all $i>1$. This relation
is clearly symmetric, and is related to the structure of the
underlying $SO_n$-representations. Indeed, when restricted to
$SO_n$, $\E{\alpha}{n} \cong \E{\beta}{n}$ if and only if $\alpha$
and $\beta$ are either associate or equal.

\vskip 10pt

\noindent{\bf ${\bf Sp_{2n}}$ Representations:} For a non-negative integer partition $\nu$ with $p$ parts where $p \leq
n$, let $\V{\nu}{2n}$ denote the irreducible rational
representation of $Sp_{2n}$ where the highest weight indexed by
the partition $\nu$ is given by the $n$ tuple:
\[
    \begin{array}{c}
        \underbrace{\left(\nu_1, \nu_2, \cdots, \nu_p, 0, \cdots, 0 \right)} \\ n
    \end{array}.
\]

\subsection{Multiplicity-Free Actions}
Let $G$ be a complex reductive algebraic group acting on a complex
vector space $V$.  We say $V$ is a {\it multiplicity-free action}
if the algebra $\calp (V)$ of polynomial functions on $V$ is
multiplicity free as a $G$ module. The criterion of
Servedio-Vinberg \cite{servedio} \cite{vinberg} says that $V$ is
multiplicity free if and only if a Borel subgroup $B$ of $G$ has a
Zariski open orbit in $V$. In other words, $B$ (and hence $G$)
acts prehomogeneously on $V$ (see \cite{sato-kimura}).  A direct
consequence is that $B$ eigenfunctions in $\calp (V)$ have a very
simple structure.  Let $Q_\psi \in \calp (V)$ be a $B$
eigenfunction with eigencharacter $\psi$, normalized so that
$Q_\psi(v_0)=1$ for some fixed $v_0$ in a Zariski open $B$ orbit
in $V$. Then $Q_\psi$ is completely determined by $\psi$: For $v =b^{-1} v_0$
in the Zariski open $B$ orbit,
$$
Q_\psi (v) = Q_\psi (b^{-1} v_0) = \psi(b) Q_\psi (v_0)= \psi (b), \qquad b\in B.
$$
$Q_{\psi}$ is then determined on all of $V$ by continuity. Since
$B = AU$, and $U = (B,B)$ is the commutator subgroup of $B$, we
can identify a character of $B$ with a character of $A$.  Thus the
$B$ eigenfunctions are precisely the $G$ highest weight vectors
(with respect to $B$) in $\calp (V)$.  Further
$$
Q_{\psi_1} Q_{\psi_2} = Q_{\psi_1 \psi_2}
$$
and so the set of $\widehat A^+ (V) = \{ \psi \in \widehat A^+
\mid Q_\psi \neq 0\}$ forms a sub-semigroup of the cone $\widehat
A^+$ of dominant weights of $A$.

An element $\psi (\neq 1)$ of a semigroup is {\it primitive} if it
is not expressible as a non-trivial product of two elements of the
semigroup. The algebra $\mcP (V)^U$ has unique factorization (see
\cite{howe-umeda}). The eigenfunctions associated to the primitive
elements of $\widehat A^+(V)$ are prime polynomials, and $\mcP
(V)^U$ is the polynomial ring on these eigenfunctions. If $\psi =
\psi_1 \psi_2$, then $Q_\psi = Q_{\psi_1} Q_{\psi_2}$. Thus, if
$\psi$ is not primitive, then the polynomial $Q_\psi$ cannot be
prime.  An element
$$
\psi = \Pi_{j=1}^k \psi_j^{c_j}
$$
has $c_j$'s uniquely determined, and hence the prime factorization
$$
Q_\psi = \Pi_{j=1}^k Q_{\psi_j}^{c_j}.
$$

Consider a multiplicity-free action of $G$ on an algebra
${\mathcal W}$.  In the general situation, we would like to
associate this algebra $\mathcal W$ with a subalgebra of
${\mathcal R} (G/U)$. With this goal in mind, we introduce the
following notion:
\begin{dfn}
Let $\mcP = \bigoplus_{\la \in \widehat A^+} \mcP_\la$ denote an
algebra graded by an abelian semigroup $\widehat A^+$. If $\mcW
\subseteq \mcP$ is a subalgebra of $\mcP$, then we say that $\mcW$
is a {\bf $(0,1)$-subalgebra of $\mcP$} if
\[
    \mcW = \bigoplus_{\la \in Z} \mcP_\la
\] where $Z$ is a sub-semigroup of $\widehat A^+$, which we will denote by $\widehat A^+(\mathcal{W}) =
Z$.  Note that $\mcW$ is graded by $\widehat A^+(\mathcal{W})$.
\end{dfn}

In what is to follow, we will usually have $\mcP = \mcP(V)$
(polynomial functions on a vector space $V$) and $\widehat A^+$ will
denote the dominant chamber of the character group of a maximal
torus $A$ of a reductive group $G$ acting on $V$.  In this
situation, we introduce an $A^+$-filtration on $\mathcal W$ as
follows:
$$
{\mathcal W}^{(\psi)} = \bigoplus_{\phi \leq \psi} {\mathcal W}_\phi
$$
where the ordering $\leq$ is the ordering on $\widehat A^+$ given by
(see \cite{popov})
$$
\begin{aligned}
\psi_1 \leq \psi_2 &\qquad \text{if $\psi_1^{-1}\psi_2$ is
expressible as a product of}\\
&\qquad \text{ rational powers of positive roots.}
\end{aligned}
$$
Note that positive roots are weights of the adjoint representation
of $G$ on its Lie algebra ${\mathfrak g}$.  We refer to the abelian
group structure on the integral weights multiplicatively. Also, it
will turn out that we only need positive \it integer \rm powers of
the positive roots.

Next consider the more specific situation where $\mathcal W$ which
is a $G$-invariant and $G$-multiplicity-free subalgebra of a
polynomial algebra $\mcP(V)$. Suppose that ${\mathcal W}^U$ has
unique factorization. Then $\mathcal W^U$ is a polynomial ring and
$\widehat A^+({\mathcal W})$ is a free sub-semigroup in $\widehat
A^+$ generated by the highest weights corresponding to the non-zero
graded components of $\mathcal W$. Write the $G$ decomposition as
follows:
$$
{\mathcal W}= \bigoplus_{\psi \in \widehat A^+({\mathcal W})} {\mathcal
W}_\psi
$$
noting that ${\mathcal W}_\psi $ is an irreducible $G$ module with
highest weight $\psi$.

If $\delta$ occurs with positive multiplicity in the tensor
product decomposition $$ {\mathcal W}_\phi \otimes  {\mathcal
W}_\psi = \bigoplus_{\delta} \dim Hom_G ({\mathcal
W}_\delta,{\mathcal W}_\phi \otimes  {\mathcal W}_\psi)\ {\mathcal
W}_{\delta},
$$
then $\delta \leq \phi \psi$.  If
$$
{\mathcal W}_\eta \subset {\mathcal W}^{(\phi)}
\quad\text{and}\quad {\mathcal W}_\gamma \subset {\mathcal
W}^{(\psi)},\quad \text{i.e., } \eta \leq \phi
\quad\text{and}\quad \gamma \leq \psi,
$$
then it follows that
$$
{\mathcal W}_\eta \cdot {\mathcal W}_\gamma \hookrightarrow
{\mathcal W}_\eta \otimes  {\mathcal W}_\gamma \subset {\mathcal
W}^{(\eta \gamma)} \subset {\mathcal W}^{(\phi \psi)}.
$$
Thus
$$
{\mathcal W}^{(\phi)}\cdot {\mathcal W}^{(\psi)} \subset {\mathcal
W}^{(\phi \psi)}.
$$
We have now an $A^+$-filtered algebra
$$
{\mathcal W} = \bigcup_{\psi \in \widehat A^+(\mathcal W) }
{\mathcal W}^{(\psi)},
$$
and this filtration is known as the {\it dominance filtration}
\cite{popov}.

With a filtered algebra, we can form its associated algebra
which is $\widehat A^+$ graded:
$$
\Gr {\mathcal W} = \bigoplus_{\psi \in \widehat A^+({\mathcal W})
} (\Gr {\mathcal W})^{\psi}
$$
where
$$
(\Gr {\mathcal W})^{\psi} = {\mathcal W}^{(\psi)} / \left(
\bigoplus_{\phi < \psi} {\mathcal W}^{(\phi)} \right).
$$

\begin{thm}
Consider a multiplicity-free $G$-module ${\mathcal W}$ with a
$\widehat A^+$-filtered algebra structure such that $\mcW$ is a
unique factorization domain.  Assume that the zero degree subspace
of ${\mathcal W}$ is $\Bbb C$.  Then there is a canonical $\widehat
A^+$-graded algebra injection:
$$
\text{$\Gr$}  \pi \ :\  \text{$\Gr$} {\mathcal W} \hookrightarrow
\calr (G/U).
$$
\end{thm}

\nt {\bf Proof.}  In \cite{howe-schur} it is shown that under the
above hypothesis, $\mcW^U$ is a polynomial ring on a canonical set
of generators. Now, $\mcW^U$ is a $\widehat A^+$-graded algebra and
therefore, there exists an injective $\widehat A^+$-graded algebra
homomorphism obtained by sending each generator of the domain to a
(indeed any) highest weight vector of the same weight in the
codomain:
$$
    \al \ :\ \mcW^U \hookrightarrow \calr (G/U)^U.
$$
Note that $\mcW^U =\Gr (\mcW^U) = (\Gr \mcW)^U$.

There exists a unique $G$-module homomorphism $\overline \al:
\Gr \mcW \hookrightarrow \calr (G/U)$ such that the following
diagram commutes:
$$
\begin{aligned}
\al:\qquad &\qquad \mcW^U &\quad \hookrightarrow &\qquad \calr (G/U)^U \\
\qquad  &\qquad\cap &\qquad \quad &\qquad \qquad \cap \\
\overline \al:\qquad  &\Gr \mcW &\qquad \hookrightarrow
&\qquad \calr (G/U)
\end{aligned}
$$

We wish to show that $\overline \al$ is an algebra homomorphism, i.e.,
$$
\begin{array}{ccccc} \\
  (\Gr \mcW)^\la & \times & (\Gr \mcW)^\mu & \genfrac{}{}{0pt}{}{\longrightarrow}{m_\mcW} & (\Gr\mcW)^{\lambda+\mu}
  \\ \\
  \overline \alpha \downarrow &  & \overline \alpha \downarrow &   & \overline \alpha \downarrow
  \\ \\
  \calr(G/U)^\lambda & \times & \calr(G/U)^\mu & \genfrac{}{}{0pt}{}{\longrightarrow}{m_{\calr(G/U)}} & \calr(G/U)^{\lambda+\mu}  \\
\end{array}
$$
commutes.

We have two maps:
$$
    f_i: (\Gr \mcW)^\lambda  \otimes  (\Gr \mcW)^\mu  \rightarrow
    \calr(G/U)^{\lambda+\mu}, \qquad i=1,2,
$$
defined by: $ f_1(v \otimes w ) = m_{\calr(G/U)}(\overline \al(v)
\otimes \overline \al(w)) $ and $ f_2(v \otimes w ) = \overline
\al (m_{\mcW}(v \otimes w)) $.

Each of $f_1$ and $f_2$ is $G$-equivariant and,
$$
\dim \left\{ \beta \ \left| \  \beta: (\Gr \mcW)^\lambda  \otimes
(\Gr \mcW)^\mu  \rightarrow \calr(G/U)^{\lambda+\mu}\right. \right\}
= 1
$$
because the Cartan product has multiplicity one in the tensor
product of two irreducible $G$-modules $V_\lambda$ and $V_\mu$
\cite{popov}.

Therefore, there exists a constant $C$ such that $f_1 = C f_2$. We
know that $\overline \al |_{\mcW^U} = \al$ is an algebra
homomorphism.  So for highest weight vectors $v^\lambda \in
\mcW^U_\lambda$ and $w^\mu \in \mcW^U_\mu$:
$$
f_1(v^\la \otimes w^\mu) = \overline \al(v^\la) \overline
\al(w^\mu) = \al(v^\la) \al(w^\mu) = \al(v^\la w^\mu) = \overline
\al(v^\la w^\mu) = f_2(v^\la \otimes w^\mu).
$$
(Note that $v^\la w^\mu$ is a highest weight vector.) Note that
$C=1$. $\qquad \square$

\subsection{Dual Pairs and Duality Correspondence}

We recall here some material and notation that will be required for
later sections.  Much of this material can be found in several other
sources on classical invariant theory (such as \cite{howe-schur} and
\cite{goodman-wallach-book}).

In our context, the theory of dual pairs may be cast in a purely
algebraic language.  In this section, we will describe three dual
pairs $(K, \fg)$, where $K$ is a classical linear algebraic group
defined over $\bbC$ and $\fg$ is a complex classical Lie algebra. In
each case, we have a linear action of $K$ on a finite dimensional
complex vector space $V$, which is a finite sum of copies of the
standard module for $K$ or copies of the dual of the standard module
for $K$. This action induces an action on the complex valued
polynomial functions on $V$, upon which $\fg$ acts by polynomial
coefficient differential operators. The actions of $\fg$ and $K$
commute with each other. Furthermore, the algebra of polynomial
coefficient differential operators which commute with the $K$-action
(resp. $\fg$-action) on $\mcP(V)$ is generated as an algebra by the
image of the $\fg$-action (resp. $K$-action). In light of this
situation, we may regard $\mcP(V)$ as a representation of $\fg$ and
$K$ simultaneously.  Theorem 3.4 describes, in part, the
(multiplicity-free) decomposition of $\mcP(V)$ into irreducible
modules for the joint action.

Of particular importance are the $K$-invariants in $\mcP(V)$.  In
each case, we may describe this invariant ring through the action of
$\fg$.  Indeed, $\fg$ may be decomposed into three subspaces denoted
$\fg^{(2,0)}$, $\fg^{(1,1)}$ and $\fg^{(0,2)}$.  These subspaces are
in fact Lie subalgebras with $\fg^{(0,2)}$ and $\fg^{(2,0)}$
abelian, while $\fg^{(1,1)}$ normalizes each of them.  Furthermore,
$\fg^{(1,1)} \oplus \fg^{(2,0)}$ and $\fg^{(1,1)} \oplus
\fg^{(0,2)}$ are the Levi decompositions of certain parabolic
subalgebras of $\fg$.

Theorem 3.3 asserts that $\mcP(V)^K$ is generated as an algebra by
$\fg^{(2,0)}$. Moreover, within a certain stable range (described in
Theorem 3.3), $\mcP(V)^K$ is isomorphic to $\mathcal S
(\fg^{(2,0)})$, the full symmetric algebra on $\fg^{(2,0)}$. This is
the \emph{first fundamental theorems of classical invariant theory}
(see \cite{goodman-wallach-book}, \cite{howe-schur}, \cite{weyl}).

Note that at the same time, we obtain an action of $K$ (by conjugation) on the
constant coefficient differential operators on $\mcP(V)$, denoted
$\mcD(V)$.  In turn, the $K$-invariant subalgebra, $\mcD(V)^K$ is
generated by $\fg^{(0,2)}$.  This brings us to our next
ingredient.  Define the $K$\emph{-harmonic polynomials} to be:
\[
    \mcH := \left\{f \in \mcP(V) \mid \Delta f = 0 \mbox{ for all } \Delta \in \mcD(V)^K \right\}.
\]
For each dual pair, we have a surjection, $\mcP(V)^K \otimes \mcH
\stackrel{m}{\rightarrow} \mcP(V)$ defined by multiplication. Note
that we may regard $\mcP(V)$ as a module over the algebra
$\mcP(V)^K$. By definition, this module is free iff $m$ is
injective. Within a certain \emph{stable range}, $m$ is indeed
injective, and this range is indicated as part of Theorem 3.3. We
have provided a proof of the injectivity part of this theorem
(also known as the Separation of Variables Theorem) as an appendix
of this paper.

For each $(\fg, K)$, the subspace $\fg^{(1,1)}$ is \emph{isomorphic}
to the Lie algebra of a subgroup $G^{(1,1)} \subseteq GL(V)$ which
commutes with the action of $K$.  In Theorem 3.5 we describe the
action of this group. Note that, in general, the differential of the
action of $G^{(1,1)}$ on $\mcP(V)$ is not quite the same as the
action of $\fg^{(1,1)}$; however, it differs only by a central
shift.

Under the joint action of $K \times G^{(1,1)}$, $\mcH$ is a
multiplicity-free invariant subspace of $\mcP(V)$.  The precise
decomposition of $\mcH$ is provided in Theorem 3.5.

Finally, in each of the three dual pair settings we have $\mcP(V)
= I(\mcJ^+) \oplus \mcH$, where $I(\mcJ^+)$ is the ideal in
$\mcP(V)$ generated by $\fg^{(2,0)}$ (which is the same as the ideal generated by the homogeneous invariants of positive
degree).  We note that the natural map:
\begin{equation}
    \mcH \longrightarrow  \mcP(V)/I(\mcJ^+)
\end{equation}
is a linear isomorphism of representations.

Details in this section including all theorems stated can be found
in \cite{goodman-wallach-book}, \cite{howe-remarks} or
\cite{howe-schur}.

\vskip 10pt

\subsubsection{Definitions of the Three Dual Pair Actions} We now describe the three dual pairs in detail as well as state
Theorems 3.3, 3.4 and 3.5 on a case-by-case basis.  For the
following, we let $M_{n,m}$ be the complex vector space of $n$ by
$m$ matrices.  We shall select a coordinate system $\{x_{ij}|
i=1,\cdots,n, \; j=1,\cdots,m\}$.

\bigskip

\nt {\bf CASE A: ${\bf (O_n, \sp_{2m})}$ where ${\bf V := M_{n, m}}$.} \\
By $O_n$ we mean the group of invertible $n \times n$ matrices, $g$
such that $g J g^t = J$ where $J$ is the $n \times n$ matrix:
\[
    J = \left(%
\begin{array}{cccc}
  0      & \cdots & 0 & 1 \\
  \vdots & \cdots & 1 & 0 \\
  0      & 1      & \cdots & \vdots \\
  1      & 0      & \cdots & 0 \\
\end{array}%
\right).
\]
This group acts on the complex $n \times m$ matrices, $V = M_{n,m}$
by left multiplication.

Using the standard matrices entries as coordinates, we define the
following differential operators:
\[
\begin{array}{lr}
\Delta_{ij} := \sum_{s=1}^n \frac{\d^2}{\d x_{si} \d x_{n-s+1 \, j}
}, \hspace{5mm} r^2_{ij} := \sum_{s=1}^n x_{si} x_{n-s+1 \, j},
\mbox{ and } & E_{ij} := \sum_{s=1}^n x_{si} \frac{\d}{\d x_{n-s+1
j}}.
\end{array}
\]
We define three spaces:
$$
\begin{aligned}
\sp_{2m}^{(1,1)} &:= \mbox{Span }\left\{ E_{ij}+ \tfrac{n}{2} \; \delta_{i,j} \mid i,j=1,\ldots, m \right\}, \\
\sp_{2m}^{(2,0)} &:= \mbox{Span }\left\{ r^2_{ij} \mid 1 \leq i \leq j \leq m \right\}, \quad \text{and} \\
\sp_{2m}^{(0,2)} &:= \mbox{Span }\left\{ \Delta_{ij} \mid 1 \leq i \leq j \leq m \right\}.
\end{aligned}
$$

The direct sum, $\fg := \sp_{2m}^{(2,0)} \oplus \sp_{2m}^{(1,1)}
\oplus \sp_{2m}^{(0,2)}$, is preserved under the usual operator
bracket and is isomorphic, as a Lie algebra, to the rank $m$ complex
symplectic Lie algebra, $\sp_{2m}$.  This presentation defines an
action of $\sp_{2m}$ on $\mcP(M_{n,m})$.  The $O_n$ action is
defined by multiplication on the left: for $g \in O_n$ and $f \in
\mcP(M_{n, m})$ we set $g \cdot f(x) = f(g^t x)$ for all $x \in
M_{n, m}$.

\bigskip

\nt {\bf CASE B: ${\bf (Sp_{2n}, \so_{2m}) }$ where ${\bf V = M_{2n, m}}$.} \\
By $Sp_{2n}$ we mean the group of complex $2n \times 2n$ invertible
matrices, $g$, such that $g J g^t = J$ where $J =
\left(%
\begin{array}{cc}
  0 & I_k \\
  -I_k & 0 \\
\end{array}%
\right)$ with $I_k$ the $k \times k$ identity matrix. This group
acts on the complex $2n \times m$ matrices, $V = M_{2n,m}$ by left
multiplication.

Using the standard matrix entries as coordinates, we define the
following differential operators:
\[
\begin{array}{c}
D_{ij} := \sum_{s=1}^n \left( \frac{\d^2}{\d x_{si} \d x_{s+n,j}} -  \frac{\d^2}{\d x_{s+n,i} \d x_{sj}}\right), \qquad S^2_{ij} := \sum_{s=1}^n (x_{si} x_{s+n,j} - x_{s+n,i} x_{sj}),
\text{ and } \\
E_{ij} := \sum_{s=1}^{2n} x_{si} \frac{\d}{\d x_{sj}}.
\end{array}
\]
We define three spaces:
$$
\begin{aligned}
\so_{2m}^{(1,1)} &:= \mbox{Span }\left\{ E_{ij}+ n \; \delta_{i,j} \mid i,j=1,\ldots, m \right\}, \\
\so_{2m}^{(2,0)} &:= \mbox{Span }\left\{ S_{ij}^2 \mid 1 \leq i < j \leq m \right\}, \qquad \text{and} \\
\so_{2m}^{(0,2)} &:= \mbox{Span }\left\{ D_{ij} \mid 1 \leq i < j \leq m \right\}.
\end{aligned}
$$

The direct sum, $\fg := \so_{2m}^{(2,0)} \oplus \so_{2m}^{(1,1)}
\oplus \so_{2m}^{(0,2)}$, is isomorphic to $\so_{2m}$, the rank $m$ orthogonal Lie algebra of type $D$, and this
presentation defines an action of $\so_{2m}$ on $\mcP(M_{2n, m})$.
For $g \in Sp_{2n}$ and $f \in \mcP(M_{2n, m})$, we set $g \cdot
f(x) = f(g^t x)$ for all $x \in M_{2n, m}$.

\bigskip

\nt {\bf CASE C: ${\bf (GL_n, \gl_{m+\ell})}$ where ${\bf V =
M_{n,m} \oplus M_{\ell,n}}$.} \\
Let $\{ x_{ab} \}$ and $\{ y_{cd}\}$ be the coordinates on
$M_{n,m}$ and $M_{\ell,n}$ respectively.  Define the following
differential operators:{\small
$$
\overline{\Delta}_{ij} := \sum_{s=1}^n \frac{\d^2}{\d x_{si} \d y_{js}}, \quad \overline{r}^2_{ij} := \sum_{s=1}^n x_{si} y_{js}, \quad E^X_{ij} := \sum_{s=1}^n x_{si} \frac{\d}{\d x_{sj}}, \quad \mbox{ and }\quad E^Y_{ij} := \sum_{s=1}^n y_{is} \frac{\d}{\d y_{js}}.
$$
} (In the above, $i$ and $j$ range over the appropriate interval
defined by the sizes of the matrices.) We define three spaces:
$$
\begin{aligned}
\gl_{m,\ell}^{(1,1)} &:= \mbox{Span }\left\{ E_{ij}^X+ \tfrac{n}{2} \; \delta_{i,j} \mid   i,j=1,\ldots, m \right\} \oplus
                     \mbox{Span }\left\{ E_{ij}^Y+ \tfrac{n}{2} \; \delta_{i,j} \mid i,j=1,\dots, \ell \right\}, \\
\gl_{m,\ell}^{(2,0)} &:= \mbox{Span }\left\{ \overline{r}^2_{ij} \mid i=1,\ldots, m, j=1,\ldots, \ell \right\}, \qquad \text{and} \\
\gl_{m,\ell}^{(0,2)} &:= \mbox{Span }\left\{ \overline{\Delta}_{ij} \mid i=1,\ldots, m,
j=1,\ldots, \ell  \right\}.
\end{aligned}
$$

The direct sum, $\fg := \gl_{m,\ell}^{(2,0)} \oplus
\gl_{m,\ell}^{(1,1)} \oplus \gl_{m,\ell}^{(0,2)}$, is isomorphic to
the rank $m+\ell$ general linear Lie algebra $\gl_{m+\ell}$, and
this presentation defines an action of $\gl_{m+\ell}$ on $\mcP(M_{n,
m} \oplus M_{\ell, n})$.  For $g \in GL_n$ and $f \in \mcP(M_{n, m}
\oplus M_{\ell, n})$, we set $g \cdot f(x,y) = f(g^t x, y
({g^t})^{-1})$ for all $x \in M_{n, m}$, $y\in M_{\ell, n}$.

\subsubsection{Theorems on the Invariants, Decompositions and Harmonics.}

Let $SM_m$ and $AM_m$ be the space of symmetric and anti-symmetric
$m$ by $m$ matrices respectively. If $V$ is a vector space, we
denote the symmetric algebra on $V$ by $\mcS(V)$. Note that in
each of the dual pairs, we have defined the action of $K$ on
$\mcP(V)$ so that $\mcP(V) \cong \mcS(V)$ as $K$ modules.  (This
is in contrast with the usual identification $\mcP(V) \cong
\mcS(V^*)$.)  Also, for a set $S$, we shall denote by $\C [ S ]$
by the algebra generated by elements in the set $S$.

\begin{thm} (First Fundamental Theorem of Invariant Theory and
Separation of Variables)
\begin{description}
\item[(a) CASE A] {\bf ${\bf (O_n, \sp_{2m})}$} The invariants
$$
\mcJ_{n,m} := \mcP(M_{ n, m})^{O_n} = \C [ r_{ij}^2] \qquad \left( \cong \calS (\sp_{2m}^{(2,0)})
\cong \mcP(SM_m) \mbox{ if } n \geq m \right).
$$
Let $\mcH_{n,m} \subseteq \mcP(M_{n,m})$ denote the
$O_n$-harmonics.  Further, if $n \geq 2m$, we have separation of variables
$$
\mcP (M_{n,m}) \simeq \mcH_{n,m} \otimes \mcJ_{n,m}.
$$
\item[(b) CASE B] {\bf ${\bf (Sp_{2n}, \so_{2m})}$} The invariants
$$
\quad \overline{\mcJ}_{2n,m} := \mcP(M_{2n, m})^{Sp_{2n}} = \C [ S_{ij}^2] \qquad  \left( \cong \calS
(\so_{2m}^{(2,0)}) \cong \mcP(AM_m) \mbox{ if } n
\geq m \right).
$$
Let $\overline{\mcH}_{2n,m} \subseteq \mcP(M_{2n,m})$ denote the
$Sp_{2n}$-harmonics.  Further, if $n \geq m$, we have separation of variables
$$
\mcP (M_{2n,m})
\simeq \overline{\mcH}_{2n,m} \otimes \overline{\mcJ}_{2n,m}.
$$
\item[(c) CASE C] {\bf ${\bf (GL_n, \gl_{m+\ell})}$} The invariants
$$
\qquad \widetilde \mcJ_{n,m,\ell} := \mcP(M_{ n, m} \oplus
M_{\ell,n})^{GL_n} = \C [ \overline{r}_{ij}^2] \qquad  \left(
\cong \calS (\gl_{m,\ell}^{(2,0)}) \cong \mcP(M_{m, \ell}) \mbox{ if } n \geq \min \, (m, \ell)
\right).
$$
Let $\widetilde{\mcH}_{n,m,\ell} \subseteq \mcP(M_{n,m} \oplus
M_{\ell,n})$ denote the $GL_n$-harmonics.  Further, if $n \geq m
+\ell$, we have separation of variables
$$
\mcP (M_{n,m} \oplus M_{\ell,n}) \simeq \widetilde
\mcH_{n, m,\ell} \otimes \widetilde \mcJ_{n, m,\ell}.
$$
\end{description}
\end{thm}

\medskip

We refer our readers to the Appendix for a new proof of the
Separation of Variables Theorem.

\medskip

The Lie algebra $\mathfrak g$ acts on $\calp (V)$ via differential
operators.  Under this action $\calp (V)$ decomposes into
irreducible (infinite-dimensional, highest weight) representations
of $\mathfrak g$.  The group $K$ is reductive, so $\mcP(V)$ also
decomposes into irreducible (finite-dimensional, highest weight)
representations of $K$.

In the space of polynomial coefficient differential operators, the
algebra generated by the image of the $K$-action is a the full
centralizer of the algebra generated by the $\mathfrak g$-action.  A
consequence of this is that $\mcP(V)$ has a multiplicity free
decomposition under the joint action of the group $K$ and the Lie
algebra $\fg$.  The irreducible constituents are of the form $U
\otimes \widetilde U$ where $U$ is an irreducible $K$-representation
while $\widetilde U$ is an irreducible $\fg$-representation.

In each of the three cases, $O_n$, $Sp_{2n}$ and $GL_n$, denote the
representations paired with $E^\la$, $V^\la$ and $F^\la$ by
$\widetilde E^\la$, $\widetilde V^\la$ and $\widetilde F^\la$
respectively.  The parametrization involving highest weights being
made precise by the pairing defined in the following theorem:

\medskip

\begin{thm} (Multiplicity-Free Decomposition under
${\bf K \times \fg}$) \newline For each case, we state the
decomposition of $\mcP(V)$ into irreducible representations:

\begin{description}
\item[(a) CASE A] {\bf ${\bf (O_n, \sp_{2m})}$}
\begin{equation}
    \mcP(M_{n, m}) = \bigoplus \E{\la}{n} \otimes \wE{\la}{2m}
\end{equation}
where the sum is over all partitions $\la$ with length at most $\min\, (n,m)$, and such that $(\la^\prime)_1+(\la^\prime)_2 \leq n$.

As a representation of $GL_m$,
\begin{eqnarray}
    \wE{\la}{2m} &= \mcJ_{n,m}  \cdot \F{\la}{m} &\qquad \mbox{for any $n, m \geq 0$,} \notag \\
                &\cong \mcS(SM_m)\otimes \F{\la}{m}       &\qquad \mbox{provided $n \geq 2m.$}
\end{eqnarray}

\item[(b) CASE B] {\bf ${\bf (Sp_{2n}, \so_{2m})}$}
\begin{equation}
    \mcP(M_{2n, m}) = \bigoplus \V{\la}{2n} \otimes \wV{\la}{2m}
\end{equation}
where the sum is over all partitions $\la$ with length at most $\min \,(n,m)$.

As a representation of $GL_m$,
\begin{eqnarray}
    \wV{\la}{2m} &= \overline{\mcJ}_{2n,m} \cdot \F{\la}{m} &\qquad \mbox{for any $n, m \geq 0$,} \notag \\
                &\cong \mcS(AM_m)\otimes \F{\la}{m}       &\qquad \mbox{provided $n \geq m.$}
\end{eqnarray}

\item[(c) CASE C] {\bf ${\bf (GL_n, \gl_{m+\ell})}$}
\begin{equation}
    \mcP(M_{n, m} \oplus M_{\ell, n}) = \bigoplus \F{(\la^+, \la^-)}{n} \otimes \wF{(\la^+,\la^-)}{m,\ell}
\end{equation}
where the sum is over all ordered pairs of partitions $(\la^+,
\la^-)$ such that \\ $\ell(\la^+) + \ell(\la^-) \leq n$,
$\ell(\la^+) \leq \min\,(n,m)$, and $\ell(\la^-) \leq
\min\,(n,\ell)$.

As a representation of $GL_m \times GL_\ell$,
\begin{eqnarray}
    \wF{(\la^+, \la^-)}{m,\ell} &= \widetilde \mcJ_{n,m,\ell} \cdot \left( \F{\la^+}{m} \otimes \F{\la^-}{\ell} \right) &\qquad \mbox{for any $n,m,\ell \geq 0$,} \notag \\
                &\cong \mcS(M_{m, \ell})\otimes \left( \F{\la^+}{m} \otimes \F{\la^-}{\ell} \right) &\qquad \mbox{provided $n \geq m+\ell.$}
\end{eqnarray}
\end{description}
\end{thm}

\vskip 10pt

\nt {\bf Remarks:} In Case C, the representation $\widetilde
F^{(\la^+, \la^-)}_{(m,\ell)}$ are (in general) complexifications of
infinite-dimensional highest weight representations of $({\mathfrak
u}_{(m,\ell)})_{\Bbb C} \simeq \gl_{m+\ell}$. Sometimes want to
emphasize the interplay of the two pieces $M_{n,m}$ and
$M_{\ell,n}$, by writing $\gl_{m,\ell}$ instead of $\gl_{m+\ell}$.
The degenerate case when $\ell=0$ is particularly interesting.
\\
This is the ${\bf GL_n \times GL_m}$ {\bf duality}:
\begin{equation}
\mcP( M_{n,m} ) =
\bigoplus_{\la} F_{(n)}^{\la} \ot F_{(m)}^{\la}
\end{equation}
where the sum is over all integer partitions $\lambda$ such that $\ell (\lambda) \leq \min\,  (n,m)$.

\vskip 10pt

\begin{thm} (Multiplicity-Free Decomposition of
Harmonics under ${\bf K \times G^{(1,1)}}$)\newline We proceed in
cases:
\begin{description}
\item[(a) CASE A] {\bf ${\bf (O_n, \sp_{2m})}\qquad $}  Let
$\mcH_{n,m} \subseteq \mcP(M_{n,m})$ denote the $O_n$-harmonics.
The group $O_n \times GL_m$ acts on $\mcP(M_{n,m})$ by $(g,h)
\cdot f(x) = f(g^t x h)$, where $g \in O_n$, $h \in GL_m$ and $x
\in M_{n,m}$.  Then $\mcH_{n,m}$ is invariant under this action.
As an $O_n \times GL_m$ representation,
\begin{equation}
\mcP(M_{n,m})/I(\mcJ_{n,m}^+) \cong    \mcH_{n,m} = \bigoplus \E{\la}{n} \otimes \F{\la}{m},
\end{equation}
where the sum is over all partitions $\la$ with length at most
$\min \,(n,m)$ and such that $(\la^\prime)_1+(\la^\prime)_2 \leq
n$.

\item[(b) CASE B] {\bf ${\bf (Sp_{2n}, \so_{2m})}\qquad $} Let $\overline{\mcH}_{2n,m} \subseteq \mcP(M_{2n,m})$ denote the $Sp_{2n}$-harmonics.\\
The group $Sp_{2n} \times GL_m$ acts on $\mcP(M_{2n,m})$ by $(g,h) \cdot f(x) = f(g^t x h)$,
where $g \in Sp_{2n}$, $h \in GL_m$ and $x \in M_{2n,m}$.  Then
$\mcH_{2n,m}$ is invariant under this action.  As a $Sp_{2n} \times GL_m$ representation,
\begin{equation}
\mcP(M_{2n,m})/I({\overline \mcJ}_{2n,m}^+) \cong    {\overline
\mcH}_{2n,m} = \bigoplus \V{\la}{2n} \otimes \F{\la}{m}
\end{equation}
where the sum is over all partitions $\la$ with length at most
$\min\,(n,m)$.

\item[(c) CASE C] {\bf ${\bf (GL_n, \gl_{m+\ell})}\qquad $} Let
$\widetilde{\mcH}_{n,m,\ell} \subseteq \mcP(M_{n,m} \oplus
M_{\ell,n})$ denote the \\ $GL_n$-harmonics.  The group $GL_n \times
GL_m \times GL_\ell$ acts on $\mcP(M_{n,m} \oplus M_{\ell,n})$ by
\[(g,h_1,h_2) \cdot f(x, y) = f(g^t x h_1, h_2^t y (g^t)^{-1}), \]
for $g\in GL_n$, $h_1 \in GL_m$, $h_2 \in GL_\ell$, $x \in M_{n,m}$
and $y \in M_{\ell, n}$. Then $\widetilde{\mcH}_{n,m,\ell}$ is
invariant under this action. As a $GL_n \times GL_m \times GL_\ell$
representation,
\begin{equation}
\mcP(M_{n,m} \oplus M_{\ell,n})/I(\widetilde{\mcJ}_{n, m,\ell}^+)
\cong    \widetilde{\mcH}_{n,m,\ell} = \bigoplus
\F{(\la^+, \la^-)}{n} \otimes \left(\F{\la^+}{m} \otimes
\F{\la^-}{\ell} \right)
\end{equation}
where the sum is over all ordered pairs of partitions $(\la^+,
\la^-)$ such that \\ $\ell(\la^+) + \ell(\la^-) \leq n$,
$\ell(\la^+) \leq \min\, (n,m)$, and $\ell(\la^-) \leq \min\, (n,
\ell)$.
\end{description}
\end{thm}

\nt {\bf Remarks.} The three cases are summarized in the following table:
\[
\begin{array}{|c||c|c|c|} \hline
\quad K \quad          & \quad O_n \quad      &  \quad Sp_{2n} \quad     & \quad GL_n \quad                  \\
\hline \fg          & \sp_{2m} & \so_{2m}  & \gl_{m,\ell}
\\ \hline V            & M_{n,m}  & M_{2n,m}  & M_{n,m}\oplus
M_{\ell,n} \\ \hline G^{(1,1)}    & GL_m     & GL_m      & GL_m
\times GL_\ell      \\ \hline
\end{array}
\]

\section{Reciprocity Algebras}

In this paper, we study branching algebras using classical invariant theory. The formulation of
classical invariant theory in terms of dual pairs \cite{howe-remarks} allows one to
realize branching algebras for classical symmetric pairs as concrete
algebras of polynomials on vector spaces. Furthermore, when realized in this
way, the branching algebras have a double interpretation in which they solve
two related branching problems simultaneously. Classical invariant theory
also provides a flexible means which allows an inductive approach to the
computation of branching algebras, and makes evident natural connections
between different branching algebras.

The easiest illustration of the above assertions is the realization of the tensor
product algebra for $GL_n$ presented as follows.

\subsection{Illustration: Tensor Product Algebra for GL$_n$}
This first example is in \cite{howe-schur}, which we recall here as
it is a model for the other (more involved) constructions of
branching algebras as total subalgebras of $GL_n$ tensor product
algebras.

Consider the joint action of $GL_n \times GL_m$ on the $\calp( M_{n,
m})$ by the rule
$$
(g, h) \cdot f(x) = f(g^t x h),\qquad \text{for } g \in GL_n, h
\in GL_m, x \in M_{n, m}.
$$
For the corresponding action on polynomials, one has the
decomposition (see Theorem 3.4(c) and (3.8))
$$
\calp (M_{n, m}) \simeq \bigoplus_{\lambda} F_{(n)}^\lambda
\otimes  F_{(m)}^\lambda , \eqno (4.1)
$$
of the polynomials into irreducible $GL_n \times GL_m$
representations. Note that the sum is over non-negative partitions
$\la$ with length at most $\min\, (n, m)$.

Let $U_m = U_{GL_m}$ denote the upper triangular unipotent
subgroup of $GL_m$. From decomposition (4.1), we can easily see
that
$$
\calp (M_{n, m})^{U_m} \simeq \left( \bigoplus_\lambda
F_{(n)}^\lambda \otimes  F_{(m)}^\lambda \right)^{U_m} \simeq
\bigoplus_\lambda F_{(n)}^\lambda \otimes  ( F_{(m)}^\lambda
)^{U_m}. \eqno (4.2)
$$
Since the spaces $(F_{(m)}^\lambda)^{U_m}$ are one-dimensional,
the sum in equation (4.2) consists of one copy of each
$F_{(n)}^\lambda$. Just as in the discussion of \S 3.2, the
algebra is graded by $\widehat A_m ^+$, where $A_m$ is the
diagonal torus of $GL_m$, and one sees from (4.2) that the graded
components are the $F_{(n)}^\lambda$.

By the arguments in \S 3.2, $\calp (M_{n,m})^{U_m}$ can thus be
associated to a graded subalgebra in $\mathcal R (GL_n/U_n)$, in
particular, this is a $(0,1)$-subalgebra as in Definition 3.1. To
study tensor products of representations of $GL_n$, we can take
the direct sum of $M_{n, m}$ and $M_{n, \ell}$. We then have an
action of $GL_n \times GL_m \times GL_{\ell}$ on $\calp (M_{n, m}
\oplus M_{n, \ell})$. Since $\calp (M_{n, m} \oplus M_{n, \ell})
\simeq \calp (M_{n, m}) \otimes  P( M_{n, \ell})$, we may deduce
from (4.1) that
$$
\calp (M_{n, m} \oplus M_{n, \ell})^{U_m \times U_{\ell}} \simeq
\calp (M_{n,m})^{U_m} \otimes  \calp (M_{n, \ell})^{U_{\ell}}
$$
$$
\simeq \bigoplus_{\mu,\nu} (F_{(n)}^\mu \otimes  F_{(n)}^\nu )
\otimes \left ( (F_{(m)}^\mu )^{U_m} \otimes  (F_{(\ell)}^\nu
)^{U_{\ell}} \right ). \eqno (4.3)
$$

Thus, this algebra is the sum of one copy of each tensor products
 $F_{(n)}^\mu \otimes  F_{(n)}^\nu $. Hence, if we take the $U_n$-invariants, we
will get a subalgebra of the tensor product algebra for $GL_n$. This results
in the algebra
$$
\left( \calp (M_{n, m} \oplus M_{n, \ell} \right)^{U_m \times U_{\ell}})^{U_n} \simeq
\calp (M_{n, m} \oplus M_{n, \ell})^{U_m \times U_{\ell} \times U_n}.
$$

This shows that we can realize the tensor product algebra for $GL_n$, or more
precisely, various $(0,1)$-subalgebras of it, as  algebras of polynomial functions
on matrices, specifically as the algebras
$\calp (M_{n, m} \oplus M_{n, \ell})^{U_m \times U_{\ell} \times U_n}$.

However, the algebra $\calp (M_{n, m} \oplus M_{n, \ell})^{U_m
\times U_{\ell} \times U_n}$ has a second interpretation, as a
different branching algebra. We note that $M_{n, m} \oplus M_{n,
\ell} \simeq M_{n, m + \ell}$.  On this space we have the action
of $GL_n \times GL_{m + \ell}$, which is described by the obvious
adaptation of equation (4.1). The action of $GL_n \times
GL_m\times GL_{\ell}$ arises by restriction of the action of
$GL_{m + \ell}$ to the  subgroup $GL_m \times GL_{\ell}$ embedded
block diagonally in $GL_{m + \ell}$. By (the obvious analog of)
decomposition (4.2), we see that
$$
\calp (M_{n, m + \ell})^{U_n} \simeq
\bigoplus_\lambda (F_{(n)}^\lambda )^{U_n} \otimes  F_{(m + \ell)}^\lambda.
$$
This algebra embeds as a subalgebra of $\calr (GL_{m + \ell} /
U_{m + \ell})$, in particular, this is a $(0,1)$-subalgebra as in
Definition 3.1. If we then take the $U_m \times U_{\ell}$
invariants, we find that
$$
(\calp (M_{n, m + \ell})^{U_n})^{U_m \times U_{\ell}}
\simeq \bigoplus_\lambda (F_{(n)}^\lambda)^{U_n} \otimes  (F_{(m + \ell)}^\lambda )^{U_m \times U_{\ell}}
$$
is (a $(0,1)$-subalgebra of) the
$(GL_{m + \ell}, GL_m \times GL_{\ell})$ branching algebra. Thus, we have
established the following result.

\begin{thm}
\begin{enumerate}

\item[(a)] The algebra $\calp (M_{n, m + \ell})^{U_n \times U_m
\times U_{\ell}}$ is isomorphic to a $(0,1)$-subalgebra of the
$(GL_n \times GL_n, GL_n)$ branching algebra (a.k.a. the $GL_n$
tensor product algebra), and to a $(0,1)$-subalgebra of the
$(GL_{m + \ell}, GL_m \times GL_{\ell})$ branching algebra.

\item[(b)] In particular,
the dimension of the $\psi^\lambda \times \psi^\mu \times \psi^\nu $ homogeneous
component for $A_n \times A_m \times A_{\ell}$ of
$\calp (M_{n, m + \ell})^{U_n \times U_m \times U_{\ell}}$ records
simultaneously
\begin{enumerate}
\item[(i)] the multiplicity of $F_{(n)}^\lambda$ in the tensor product
$F_{(n)}^\mu \otimes  F_{(n)}^\nu $, and
\item[(ii)] the multiplicity of
$F_{(m)}^\mu \otimes  F_{(\ell)}^\nu $ in $F_{(m + \ell)}^\lambda$,
\end{enumerate}
for partitions $\mu$, $\nu$, $\lambda$ such that $\ell (\mu) \leq
\min(n,m)$, $\ell(\nu) \leq \min (n, \ell)$ and $\ell (\lambda)
\leq \min (n, m+\ell)$.

\end{enumerate}
\end{thm}

Thus, we can not only realize the $GL_n$ tensor product algebra concretely as an
algebra of polynomials, we find that it appears simultaneously in two
guises, the second being as the branching algebra for the pair
$(GL_{m + \ell}, GL_m \times GL_{\ell})$. We emphasize two features of this
situation.

First, the pair $(GL_{m + \ell}, GL_m \times GL_{\ell} )$, as well as the
pair $(GL_n \times GL_n, GL_n)$, is a symmetric pair. Hence, both the
interpretations of
$\calp (M_{n, m + \ell})^{U_n \times U_m \times U_{\ell}}$ are as branching
algebras for symmetric pairs.

Second, the relationship between the two situations is captured by
the notion of \lq \lq see-saw pair" of dual pairs
\cite{kudla-see-saw}. Precisely, a context for understanding the
decomposition law (4.1) is provided by observing that $GL_n$ and
$GL_m$ (or more correctly, slight modifications of their Lie
algebras) are mutual centralizers inside the Lie algebra $\spm
(M_{n, m})$ (of the metaplectic group) of polynomial coefficient
differential operators of total degree two on $M_{n, m}$
\cite{howe-remarks} \cite{howe-schur}. We say that they define a
\it dual pair \rm inside $\spm (M_{n, m})$. The decomposition
(4.1) then appears as the correspondence of representations
associated to this dual pair \cite{howe-remarks}. Further, the
pairs of groups $(GL_n, GL_{m + \ell}) = (G_1, G_1')$ and $(GL_n
\times GL_n, GL_m \times GL_{\ell}) = (G_2, G_2')$ both define
dual pairs inside the Lie algebra $\spm (M_{n, m + \ell})$. We
evidently have the relations
$$
G_1 = GL_n \subset GL_n \times GL_n = G_2, \eqno (4.4)
$$
and (hence)
$$
G_1' = GL_ {m + \ell} \supset GL_m \times GL_{\ell} = G_2'. \eqno
(4.5)
$$
We refer to a pair of dual pairs related as in inclusions (4.4)
and (4.5), a \it see-saw pair \rm of dual pairs.

In these terms, we may think of the symmetric pairs $(G_2, G_1)$ and
$(G_1', G_2')$ as a ``{\it reciprocal pair}'' of symmetric pairs. If
we do so, we see that the algebra $\calp (M_{n, m + \ell})^{U_n
\times U_m \times U_{\ell}}$ is describable as $\calp (M_{n, m +
\ell})^{U_{G_1} \times U_{G'_2}}$ -- it has a description in terms
of the see-saw pair, and in this description the two pairs of the
see-saw, or alternatively, the two reciprocal symmetric pairs, enter
equivalently into the description of the algebra that describes the
branching law for both symmetric pairs. For this reason, we also
call this algebra, which describes the branching law for both
symmetric pairs, the {\it reciprocity algebra} of the pair of pairs.

It turns out that any branching algebra associated to a classical symmetric
pair, that is, a pair $(G, H)$ in which $G$ is a product of classical
groups, has an interpretation as a reciprocity algebra -- an algebra that
describes a branching law for two reciprocal symmetric pairs simultaneously.
Sometimes, however, one of the branching laws involves infinite-dimensional representations.

\subsection{Symmetric Pairs and Reciprocity Pairs}
In the context of dual pairs, we would like to understand the $(G,
H)$ branching of irreducible representations of $G$ to $H$, for
symmetric pairs $(G, H)$. Table I lists the symmetric pairs which
we will cover in this paper.

If $G$ is a classical  group over $\C$, then $G$ can be embedded as
one member of a dual pair in the symplectic group as described in
\cite{howe-remarks}. The resulting pairs of groups are $(GL_n , GL_m
)$ or $(O_n , Sp_{2m})$, each inside $Sp_{2nm}$, and are called {\it
irreducible} dual pairs. In general, a dual pair of reductive groups
in $Sp_{2r}$ is a product of such pairs.

\newpage
\begin{center} {\bf Table I: Classical Symmetric Pairs} \end{center}
\vskip 10pt
\begin{center}
\begin{tabular}{|c|c|c|} \hline
{\qquad \bf Description\qquad}  & $\qquad {\bf G}\qquad $ & $\qquad {\bf H}\qquad $ \\ \hline
Diagonal  &\qquad $GL_n \times GL_n\qquad$ & $GL_n$ \\ \hline
Diagonal  & $O_n \times O_n$ & $O_n$ \\ \hline
Diagonal & $Sp_{2n} \times Sp_{2n} $ & $Sp_{2n}$ \\ \hline
Direct Sum & $GL_{n+m}$ &\qquad $GL_n \times GL_m\qquad$ \\ \hline
Direct Sum & $O_{n+m}$ & $O_n \times O_m$ \\ \hline
Direct Sum & $Sp_{2(n+m)}$ & $Sp_{2n} \times Sp_{2m}$ \\ \hline
Polarization & $O_{2n}$ & $GL_n$ \\ \hline
Polarization & $Sp_{2n}$ & $GL_n$ \\ \hline
Bilinear Form & $GL_n$ & $O_n$  \\ \hline
Bilinear Form & $GL_{2n}$ & $Sp_{2n}$  \\ \hline
\end{tabular}
\end{center}
\vskip 10pt

\vskip 10pt

\nt {\bf Proposition 4.2} {\it Let $G$ be a classical group, or a
product of two copies of a classical group. Let $G$ belong to a
dual pair $(G, G')$ in a symplectic group $Sp_{2m}$.  Let $H
\subset G$ be a symmetric subgroup, and let $H'$ be the
centralizer of $H$ in $Sp_{2m}$. Then $(H, H')$ is also a dual
pair in $Sp_{2m}$, and $G'$ is a symmetric subgroup inside $H'$. }

\vskip 10pt
\nt {\bf Proof:} This can be shown by fairly easy case-by-case checking. The basic
reason that $(H, H')$ form a dual pair is that, for any classical symmetric
pair $(G, H)$, the restriction of the standard module of
$G$, or its dual, to $H$ is a sum of standard modules of $H$, or their
duals \cite{howe-remarks}. This is very easy to check on a case-by-case basis. The see-saw
relationship of symmetric pairs organizes the 10 series of symmetric pairs
as given in Table I into five pairs of pairs. These are shown in Table II. $\qquad \square$

\begin{center} {\bf Table II: Reciprocity Pairs}\end{center}

\medskip
{\small
\begin{center}
\begin{tabular}{|c|c|c|} \hline
\ {\bf Symmetric Pair } ${\bf (G,H)}$ \quad &\quad ${\bf
(H,{\mathfrak h}')}$ \quad &\quad ${\bf (G,{\mathfrak g}')}$ \quad  \\
\hline \hline \ $(GL_n \times GL_n, GL_n)$ \ &\ $(GL_n,
\gl_{m+\ell})$ \ &\ $(GL_n \times GL_n, \gl_m \oplus \gl_l)$ \ \\
\hline \ $(O_n \times O_n, O_n)$ \ &\ $(O_n, \spm_{2(m+\ell )})$ \
&\ $(O_n \times O_n, \spm_{2m} \oplus \spm_{2l})$ \ \\ \hline \
$(Sp_{2n} \times Sp_{2n}, Sp_{2n})$ \ &\ $(Sp_{2n},
\so_{2(m+\ell)})$ \ &\
$(Sp_{2n} \times Sp_{2n}, \so_{2m} \oplus  \so_{2\ell})$ \ \\
\hline \ $(GL_{n+m}, GL_n \times GL_m)$ \ &\ $(GL_n \times GL_m,
\gl_\ell \oplus \gl_\ell)$ \ &\ $(GL_{n+m}, \gl_\ell)$ \ \\
\hline \ $(O_{n+m}, O_n \times O_m)$ \ &\ $(O_n \times O_m,
\spm_{2\ell} \oplus \spm_{2\ell})$ \ &\ $(O_{n+m}, \spm_{2\ell})$
\ \\  \hline \ $(Sp_{2(n+m)}, Sp_{2n} \times Sp_{2m})$ \ &\
$(Sp_{2n} \times Sp_{2m}, \so_{2\ell } \oplus \so_{2\ell })$  \ &\
$(Sp_{2(n+m)}, \so_{2\ell })$ \ \\ \hline \ $(O_{2n}, GL_n)$ \ &\
$(GL_n, \gl_{m,m})$ \ &\ $(O_{2n}, \spm_{2m})$ \ \\ \hline \
$(Sp_{2n},
GL_n)$ \ &\ $(GL_n, \gl_{m,m})$ \ &\ $(Sp_{2n}, \so_{2m})$ \ \\
\hline \quad $(GL_{n}, O_n)$ \ &\ $(O_n, \spm_{2m})$ \ &\
$(GL_{n}, \gl_m)$ \ \\ \hline \ $(GL_{2n}, Sp_{2n})$ \ &\
$(Sp_{2n}, \so_{2m})$ \ &\ $(GL_{2n}, \gl_{m})$ \quad \\ \hline
\end{tabular}
\end{center}
}
\medskip
\noindent{\bf Remark:} Note that when the second component of any
pair in Table II is of Lie type A, then the action actually
integrates to the group.  Table II also amounts to another point of
view on the structure on which \cite{howe-transcending} is based.

\vskip 10pt

We begin with discussions of reciprocity algebras in the next
three sections.  The discussions provided are ordered more in
terms of complexity and do not follow the sequence given in Table
I.

\section{Branching from $GL_n$ to $O_n$}

Consider the problem of restricting irreducible representations of
$GL_n$ to the orthogonal group $O_n$. We consider the symmetric
see-saw pair $(GL_n, O_n)$ and $(Sp_{2m}, GL_m)$. As in the
discussion of \S 4.1, we can realize (a $(0,1)$-subalgebra of) the
coordinate ring of the flag manifold $GL_n/U_n$ as the algebra of
$U_m$-invariants on $\calp (M_{n, m})$. If we then look at the
$U_{O_n}$-invariants in this algebra, then we will have (a certain
$(0,1)$-subalgebra of) the $(GL_n, O_n)$ branching algebra. Thus, we
are interested in the algebra
$$
\calp (M_{n,m})^{U_{O_n} \times U_m}.
$$
We note that, in analogy with the situation of \S 4.1, this is the
algebra of invariants for the unipotent subgroups of the smaller
member of each symmetric pair.

Let us investigate what this algebra appears to be if we first take
invariants with respect to $U_{O_n}$.  We have a decomposition of
$\calp (M_{n, m})$ as a joint $O_n \times \spm_{2m}$-module (see
Theorem 3.4 (a)):
$$
\calp (M_{n, m}) \simeq \bigoplus_\mu E_{(n)}^\mu \otimes
\widetilde E_{(2m)}^{\mu}. \eqno (5.1)
$$
Recall that the sum runs through the set of all non-negative
integer partitions $\mu$ such that $\ell (\mu) \leq \min(n,m)$ and
$(\mu')_1 + (\mu')_2 \leq n$. Here $E_{(n)}^\mu$ denotes the
irreducible $O_n$ representation parameterized by $\mu$.  Recall
from \S 3.3, the decomposition $\calp (M_{n, m}) \simeq
\bigoplus_\mu F_{(n)}^\mu \otimes F_{(m)}^{\mu}$.  The module
$E_{(n)}^\mu$ is generated by the $GL_n$ highest weight vector in
$F_{(n)}^\mu$. Further, $\widetilde E_{(2m)}^{\mu}$ is an
irreducible infinite-dimensional representation of $\spm_{2m}$
with lowest $\gl_m$-type $F_{(m)}^\mu$.

\begin{thm} Assume $n > 2m$.
\begin{enumerate}
\item[(a)] The algebra $\calp (M_{n, m})^{U_{O_n} \times U_m }$ is
isomorphic to a $(0,1)$-subalgebra of the  $(GL_n, O_n)$ branching
algebra, and to a $(0,1)$-subalgebra of the $(\sp_{2m}, GL_m)$
branching algebra.

\item[(b)] In particular,
the dimension of the $\phi^\mu \times \psi^\la$ homogeneous
component for $A_{O_n} \times A_m$ of
$\calp (M_{n, m})^{U_{O_n} \times U_m}$ records
simultaneously
\begin{enumerate}
\item[(i)] the multiplicity of $E_{(n)}^\mu$ in the representation
$F_{(n)}^\la$, and
\item[(ii)] the multiplicity of $F_{(m)}^\lambda$ in $\widetilde E_{(2m)}^\mu$.
\end{enumerate}
for partitions $\mu$, $\lambda$ such that $\ell (\mu) \leq m$, $\ell
(\lambda) \leq m$.

\end{enumerate}
\end{thm}

\nt {\bf Proof.} Taking the $U_{O_n}$-invariants for the
decomposition (5.1), we find that
$$
\calp (M_{n, m})^{U_{O_n}} \simeq \bigoplus_\mu
(E_{(n)}^\mu)^{U_{O_n}} \otimes \widetilde E_{(2m)}^{\mu}, \eqno
(5.2)
$$
where the sum is over partitions $\mu$ such that $\ell (\mu) \leq
\min (n,m)$ and $(\mu')_1 + (\mu')_2 \leq n$.  Note that the
stability condition $n>2m$ guarantees the latter inequality.  The
space $(E_{(n)}^\mu)^{U_{O_n}}$ is the space of highest weight
vectors for $(E_{(n)}^\mu)^{U_{O_n}}$. We would like to say that it
is one-dimensional, so that $\calp (M_{n, m})^{U_{O_n}}$ would
consist of one copy of each of the irreducible representations
$\widetilde E^{\mu}_{(2m)}$. But, owing to the disconnectedness of
$O_n$, this is not quite true, and when it is true, the highest
weight may not completely determine $E_{(n)}^\mu$.

However, if $n > 2m$, then $(E_{(n)}^\mu)^{U_{O_n}}$ is
one-dimensional, and does single out $E_{(n)}^\mu$ among the
representations which appear in the sum (5.1). Hence, let us make
this restriction for the present discussion. Taking the $U_m$
invariants in the sum (5.2), we find that
$$
(\calp (M_{n, m})^{U_{O_n}})^{U_m} \simeq \bigoplus_\mu
(E_{(n)}^\mu)^{U_{O_n}} \otimes (\widetilde E^{\mu}_{(2m)})^{U_m}.
\eqno (5.3)
$$
Note that the sum is over all partitions $\mu$ such that $\ell (\mu)
\leq m$ (since $n > 2m$).  The space $(\widetilde
E^{\mu}_{(2m)})^{U_m}$ describes how the representation $\widetilde
E^{\mu}_{(2m)}$ of $\spm_{2m}$ decomposes as a $\gl_m$ module, or
equivalently, as a $GL_m$-module. In other words, $(\widetilde
E^{\mu}_{(2m)})^{U_m}$ describes the branching rule from $\spm_{2m}$
to $\gl_m$ for the module $\widetilde E^{\mu}_{(2m)}$.

We know (thanks to our restriction to $n > 2m$) that the space
$(E_{(n)}^\mu)^{U_{O_n}}$ is one-dimensional. Let $\phi^\mu$ be
the $A_{O_n}$ weight of $( E_{(n)}^\mu)^{U_{O_n}}$. Thus,
$\phi^\mu$ is the restriction to the diagonal maximal torus
$A_{O_n}$ of the character $\psi^\mu$ of the group $A_n$ of
diagonal $n \times n$ matrices. Our assumption further implies
that $\phi^\mu$ determines $E_{(n)}^\mu$. Therefore, for a given
dominant $A_m$ weight $\psi^\lambda$, corresponding to the
partition $\lambda$, where $\ell (\lambda) \leq m$, the
$\psi^\lambda $-eigenspace in $(\widetilde E^{\mu}_{(2m)})^{U_m}$
tells us the multiplicity of $F_{(m)}^\lambda $ in the restriction
of $\widetilde E^{\mu}_{(2m)}$ to $\gl_m$. This is the same as the
dimension of the joint ($\phi^\mu \times \psi^\lambda$)-eigenspace
in
$$
(\calp (M_{n, m})^{U_{O_n}})^{U_m}\simeq \calp (M_{n, m})^{U_{O_n} \times U_m} \simeq
(\calp (M_{n, m})^{U_m})^{U_{O_n}}.
$$

But we have already seen that this
eigenspace describes the multiplicity of $E_{(n)}^\mu$ in $F_{(n)}^\lambda$. Thus,
again the $A_{O_n} \times A_m$ homogeneous components of
$\calp (M_{n, m})^{U_{O_n} \times U_m}$ have a simultaneous interpretation, one
for a branching law associated to each of the two symmetric pairs composing
the symmetric see-saw pair. $\qquad \qquad \square$

In this case, one of the branching laws
involves infinite-dimensional representations. However, they are highest weight representations, which are the most tractable of
infinite-dimensional representations, from an algebraic point of view.

\section{Tensor Product Algebra for $O_n$}

Using the symmetric see-saw pair $\left( (O_n, O_n \times O_n),
(Sp_{2m} \times Sp_{2 \ell}, Sp_{2(m + \ell)}) \right)$, we can
construct ($(0,1)$-subalgebras of) the tensor product algebra for
$O_n$. To prepare for this, we should explicate the decomposition
(5.1) further.

Let us recall the basic setup as in \S 3.3.1 Case A.  Recall that
${\mathcal J}_{n,m} = \calp (M_{n,m})^{O_n}$ is the algebra of
$O_n$-invariant polynomials. Theorem 3.3(a) implies that
${\mathcal J}_{n,m}$ is a quotient of $\calS (\spm_{2m}^{(2,0)})$,
the symmetric algebra on $\spm_{2m}^{(2,0)}$.

The natural mapping
$$
\Har_{n,m} \rightarrow \calp (M_{n,m})/I({\mathcal J}_{n,m}^+)
\eqno (6.1)
$$
is a linear $O_n \times GL_m$-module isomorphism.  Further, the
$O_n \times GL_m$ structure of $\Har_{n,m}$ is as follows (see
Theorem 3.5(a)):
$$
\Har_{n,m} \simeq \bigoplus_\mu E_{(n)}^\mu \otimes  F_{(m)}^\mu.
\eqno (6.2)
$$
Here $\mu$ ranges over the same diagrams as in (5.1).

From Theorem 3.4(a),
$$
\widetilde E_{(2m)}^{\mu} \simeq F_{(m)}^\mu \cdot {\mathcal
J}_{n,m} \simeq \calS (\sp_{2m}^{(2,0)}) \cdot F_{(m)}^\mu, \eqno
(6.3)
$$
and it follows that
$$
\widetilde E_{(2m)}^{\mu}/(\spm_{2m}^{(2,0)} \cdot \widetilde E_{(2m)}^{\mu} )
\simeq F_{(m)}^\mu.
$$
In other words, we can detect the $\spm_{2m}$ isomorphism class of the module
$\widetilde E_{(2m)}^{\mu}$ by the $GL_m$ isomorphism class of the quotient
$\widetilde E_{(2m)}^{\mu}/(\spm_{2m}^{(2,0)} \cdot \widetilde E_{(2m)}^{\mu} )$.
Also, if $W \subset \calp (M_{n,m})$ is any
$\spm_{2m}$-invariant subspace, then
$$
W/(\spm_{2m}^{(2,0)} \cdot W) \simeq W \cap \Har_{n,m},
$$
and this subspace also
reveals the $\spm_{2m}$ isomorphism type of $W$.

We can use the above to find a model for (a $(0,1)$-subalgebra of) the tensor product algebra of $O_n$. One consequence of the above discussion is that
$$
\left ( \calp (M_{n,m})/I({\mathcal J}_{n,m}^+) \right )^{U_m}
\simeq \bigoplus_\mu  E_{(n)}^\mu \otimes (F^\mu_{(m)})^{U_m}
\eqno (6.4)
$$
consists of one copy of each irreducible representation $E_{(n)}^\mu$.

If we repeat the above discussion for $M_{n, \ell}$, and combine the
results, we find that
$$
\left ( \calp (M_{n,m})/I({\mathcal J}_{n,m}^+) \right )^{U_m} \otimes
\left ( \calp (M_{n,\ell})/I({\mathcal J}_{n,\ell}^+) \right )^{U_{\ell}}
$$
$$
\simeq \bigoplus_{\mu, \nu} \left ( E_{(n)}^\mu \otimes
E_{(n)}^\nu \right ) \otimes  \left ( (F_{(m)}^\mu)^{U_m} \otimes
(F_{(\ell)}^\nu)^{U_{\ell}} \right ) \eqno (6.5)
$$
is a direct sum of one copy of each possible tensor product of an
$E_{(n)}^\mu$ with an $E_{(n)}^\nu$.  At this point, we make the
assumption that $n > 2(m+\ell)$, as in this range the $O_n$
constituents of decomposition are irreducible when restricted to the
connected component of the identity in $O_n$. If we now take the
$U_{O_n}$-invariants in equation (6.5), we will have (a
$(0,1)$-subalgebra of) the tensor product algebra of $O_n$:
$$
\left ( (\calp (M_{n,m})/I({\mathcal J}_{n,m}^+) )^{U_m} \otimes
 ( \calp (M_{n,\ell})/I({\mathcal J}_{n,\ell}^+)) ^{U_{\ell}}\right )^{U_{O_n}}
$$
$$
\simeq \bigoplus_{\mu, \nu} \left ( E_{(n)}^\mu \otimes
E_{(n)}^\nu \right )^{U_{O_n}} \otimes  \left (
(F_{(m)}^\mu)^{U_m} \otimes  (F_{(\ell)}^\nu)^{U_{\ell}} \right ).
\eqno (6.6)
$$

We can describe this algebra in another way. Begin with the observation that
$\calp (M_{n,m}) \otimes  \calp (M_{n, \ell}) \simeq \calp (M_{n, m + \ell})$, and
$$
\calp (M_{n,m})/I({\mathcal J}_{n,m}^+) \otimes  \calp (M_{n, \ell})/I({\mathcal J}_{n,\ell}^+)
\simeq \calp (M_{n, m + \ell})/I({\mathcal J}_{n,m}^+ \oplus {\mathcal J}_{n,\ell}^+).
$$
Thus
$$
(\calp (M_{n,m})/I({\mathcal J}_{n,m}^+))^{U_m} \otimes  (\calp (M_{n, \ell})/I({\mathcal J}_{n,\ell}^+))^{U_{\ell}}
\simeq (\calp (M_{n, m + \ell})/I({\mathcal J}_{n,m}^+ \oplus {\mathcal J}_{n,\ell}^+))^{U_m \times U_{\ell}},
$$
and
$$
\left ( (\calp (M_{n,m})/I({\mathcal J}_{n,m}^+))^{U_m} \otimes
(\calp (M_{n, \ell})/I({\mathcal J}_{n,\ell}^+))^{U_{\ell}} \right )^{U_{O_n}}
$$
$$
\simeq \left ( (\calp (M_{n, m + \ell})/I({\mathcal J}_{n,m}^+ \oplus {\mathcal J}_{n,\ell}^+))^{U_m \times U_{\ell}}\right )^{U_{O_n}}
$$
$$
\simeq
\left ( (\calp (M_{n, m + \ell})/I({\mathcal J}_{n,m}^+ \oplus {\mathcal J}_{n,\ell}^+))^{U_{O_n}} \right )^{U_m \times U_{\ell}}.
$$
\newpage
\begin{thm}  Given positive integers $n$, $m$ and $\ell$ with $n
> 2(m+\ell)$ we have:
\begin{enumerate}
\item[(a)] The algebra
$$
\left ( (\calp (M_{n,m})/I({\mathcal J}_{n,m}^+) )^{U_m} \otimes
 ( \calp (M_{n,\ell})/I({\mathcal J}_{n,\ell}^+)) ^{U_{\ell}}\right )^{U_{O_n}} $$
$$
\simeq \left ((\calp (M_{n,m + \ell})/I({\mathcal J}_{n,m}^+ \oplus {\mathcal J}_{n,\ell}^+) )^{U_{O_n}} \right ) ^{U_m \times U_{\ell}}
$$
is isomorphic to a $(0,1)$-subalgebra of the $(O_n \times O_n, O_n)$
branching algebra (a.k.a. the $O_n$ tensor product algebra), and to
a $(0,1)$-subalgebra of the $(\sp_{2(m + \ell)}, \sp_{2m} \oplus
\sp_{2\ell})$ branching algebra.

\item[(b)] Specifically, the dimension of the ($\phi^\lambda \times \psi^\mu \times \psi^\nu$)-eigenspace for $A_{O_n} \times A_m \times A_{\ell}$ of
$\left( (\calp (M_{n,m + \ell})/I({\mathcal J}_{n,m}^+ \oplus {\mathcal J}_{n,\ell}^+) )^{U_{O_n}} \right) ^{U_m \times U_{\ell}}$ records simultaneously
\begin{enumerate}
\item[(i)] the multiplicity of
$E_{(n)}^\lambda$ in $E_{(n)}^\mu \otimes  E_{(n)}^\nu$, as well as
\item[(ii)] the multiplicity of $\widetilde E_{(2m)}^{\mu} \otimes   \widetilde E_{(2\ell)}^{\nu}$ in the restriction of $\widetilde E_{(2(m+ \ell))}^{\lambda}$.
\end{enumerate}
Here the partitions $\mu$, $\nu$, $\lambda$ satisfy the following
conditions: \\
$\ell (\mu) \leq \min (n,m)$, $\ell (\nu) \leq \min (n,\ell)$, and
$\ell (\lambda) \leq \min (n,m+\ell)$.
\end{enumerate}
\end{thm}

\nt {\bf Proof.} Let us now compute the ring expressed in this
way. From Theorem 3.4(a), we know that
$$
\calp (M_{n,m})^{U_{O_n}} \simeq \left(  \bigoplus_\mu E_{(n)}^\mu \otimes  \widetilde E_{(2m)}^\mu \right)^{U_{O_n}} \simeq \bigoplus_\mu ( E_{(n)}^\mu )^{U_{O_n}} \otimes   \widetilde E_{(2m)}^\mu .
$$
Note that within the range $n > 2(m + \ell)$ we have $\dim (
E_{(n)}^\mu )^{U_{O_n}} = 1$ since the $O_n$-representations
$E_{(n)}^\mu$ remain irreducible when restricted to $SO_n$.

Now repeat this with $m$ replaced by $m + \ell$:
$$
\calp (M_{n,m+\ell})^{U_{O_n}} \simeq \left(  \bigoplus_\mu E_{(n)}^\mu \otimes  \widetilde E_{(2(m+\ell))}^\mu \right)^{U_{O_n}} \simeq \bigoplus_\mu ( E_{(n)}^\mu )^{U_{O_n}} \otimes   \widetilde E_{(2(m+\ell))}^\mu .
$$
Hence
$$
(\calp (M_{n,m + \ell})/I({\mathcal J}_{n,m}^+ \oplus {\mathcal J}_{n,\ell}^+))^{U_{O_n}} \simeq
\left( \left( \bigoplus_\lambda E_{(n)}^\lambda \otimes  \widetilde E_{(2(m + \ell))}^\lambda \right) /I({\mathcal J}_{n,m}^+ \oplus {\mathcal J}_{n,\ell}^+) \right)^{U_{O_n}}
$$
$$
\simeq  \bigoplus_\lambda ( E_{(n)}^\lambda)^{U_{O_n}}
\otimes  \left ( \widetilde E_{(2(m + \ell))}^{\lambda} / (\spm_{2m}^{(2,0)} \oplus
\spm_{2 \ell}^{(2,0)})\cdot \widetilde E_{(2(m + \ell))}^{\lambda} \right).
$$
From this we finally get
$$
\left ((\calp (M_{n,m + \ell})/I({\mathcal J}_{n,m}^+ \oplus
{\mathcal J}_{n,\ell}^+) )^{U_{O_n}} \right)^{U_m \times U_{\ell}}
$$
$$
\simeq  \bigoplus_\lambda( E_{(n)}^\lambda )^{U_{O_n}}
\otimes  \left ( \widetilde E_{(2(m + \ell))}^{\lambda} /(\spm_{2m}^{(2,0)} \oplus \spm_{2 \ell}^{(2,0)}) \cdot \widetilde E_{(2(m + \ell))}^{\lambda}  \right)^{U_m \times U_{\ell}}.
$$

From the discussion following equation (6.1), we see that the
factor
$$
\left ( \widetilde E_{(2(m + \ell))}^{\lambda } /(\spm_{2m}^{(2,0)} \oplus
\spm_{2 \ell}^{(2,0)})\cdot \widetilde E_{(2(m + \ell))}^{\lambda} \right)^{U_m \times U_{\ell}}
$$
tells us the
$\spm_{2m} \oplus \spm_{2 \ell}$ decomposition of $\widetilde E_{(2(m + \ell))}^{\lambda}$. $\qquad \square$

Hence, again the algebra has a double interpretation, one in terms of decomposing
tensor products of $O_n$ representations, and one in terms of branching from
$\spm_{2(m + \ell)}$ to $\spm_{2m } \oplus \spm_{2 \ell}$ (although the second
branching law involves infinite-dimensional representations).

\section{More Reciprocity Algebras for $(GL_n, GL_m)$}

Whereas our first example of a reciprocity algebra in \S 4.1
involved only finite-dimensional representations, the others all
involve infinite-dimensional representations in some respect. It
turns out that the apparently exceptional nature of the
reciprocity algebra for the pair $(GL_n, GL_m)$ is somewhat
deceptive. In fact, we can associate several reciprocity algebras
to $(GL_n, GL_m)$, and nearly all of them will involve
infinite-dimensional representations.
\smallskip

We shall refer to \S 3.3.1 and consider the action of $GL_n$ on
$\calp (M_{n, m} \oplus M_{\ell, n})$ by the rule
$$
g \cdot f(x, y) = f(g^t x, y(g^t)^{-1}) \eqno (7.1)
$$
for $x \in M_{n,m}$, $y \in M_{\ell,n}$ and $g \in GL_n$. Recall
from Theorem 3.3(c) that the algebra $\widetilde{\mathcal
J}_{n,m,\ell}$ generated by $\gl_{m , \ell}^{(2,0)}$. It is the
space of all polynomials on $\calp (M_{n, m} \oplus M_{\ell, n})$
invariant under $GL_n$. Let $\widetilde \Har_{n,m,\ell} $ be the
space of $GL_n$-harmonics and recall the $GL_n \times GL_m \times
GL_\ell$ isomorphism (see Theorem 3.5(c)):
$$
\calp (M_{n, m} \oplus M_{\ell, n})/I(\widetilde{\mathcal J}_{n,m, \ell}^+ )\simeq \widetilde  \Har_{n,m,\ell}.
$$
\begin{thm}
\begin{enumerate}
\item[(a)] The algebra
$$
\left ( ( \calp (M_{n, m} \oplus M_{\ell, n})/I(\widetilde{\mathcal J}_{n,m, \ell}^+ ) )^{U_m \times U_{\ell}} \otimes
( \calp (M_{n, m'} \oplus M_{\ell', n})/I(\widetilde{\mathcal J}_{n,m',\ell'}^+))^{U_{m'} \times U_{\ell'}}
\right ) ^{U_n}
$$
is isomorphic to a $(0,1)$-subalgebra of the $(GL_n\times
GL_n,GL_n)$ branching algebra as well as to a
$(\gl_{m+m',\ell+\ell'}, \gl_{m,\ell} \oplus \gl_{m',\ell'})$
branching algebra.

\item[(b)] In particular,  the dimension of the ($A_n \times A_m \times A_{m'} \times A_{\ell} \times A_{\ell'}$)-eigenspace of\linebreak $\calp (M_{n, m+m'} \oplus M_{\ell + \ell', n})$ describes simultaneously
\begin{enumerate}
\item[(i)] the multiplicity of $F_{(n)}^{(\lambda^+, \lambda^-)}$ in
$ F_{(n)}^{(\mu^+ , \mu^-)} \otimes  F_{(n)}^{(\nu^+, \nu^-)}$, and
\item[(ii)] the multiplicity of the representation
$\widetilde F_{(m, \ell)}^{(\mu^+,\mu^-)} \otimes  \widetilde F_{(m', \ell')}^{(\nu^+, \nu^-)}$ of $\gl_{m ,\ell} \oplus \gl_{m',\ell'}$ in the restriction of the
representation
$\widetilde F_{(m+m', \ell+\ell')}^{(\lambda^+, \lambda^-)}$ of $\gl_{(m+m'), (\ell + \ell')}$.
\end{enumerate}
Here the partitions $\mu^+$, $\mu^-$, $\nu^+$, $\nu^-$, $\la^+$,
$\la^-$ are such that $\ell (\mu^+) \leq m$, $\ell (\mu^-) \leq \ell
$, $\ell (\mu^+)+\ell (\mu^-) \leq n$, $\ell (\nu^+) \leq m'$, $\ell
(\nu^-) \leq \ell'$, $\ell (\nu^+)+ \ell (\nu^-) \leq n$, $\ell
(\la^+) \leq m+m'$, $\ell (\la^-) \leq \ell + \ell'$ and $\ell
(\la^+)+ \ell (\la^-) \leq n$.

\end{enumerate}
\end{thm}

\vskip 5pt

\nt {\bf Remarks:} Recall from the remarks after Theorem 3.4 that
we have written $\gl_{m,\ell}$ instead of $\gl_{m+\ell}$ to
emphasize the interplay of the two components $M_{n,m}$ and
$M_{\ell,n}$.

\vskip 5pt

\nt {\bf Proof.} From the above description of $\calp (M_{n, m}
\oplus M_{\ell, n})$, we can see using (3.11) that
$$
( \calp (M_{n, m} \oplus M_{\ell, n})/I(\widetilde{\mathcal J}_{n,m, \ell}^+ ) )^{U_m \times U_{\ell}}
\simeq  \bigoplus_{\mu^+, \mu^-} F_{(n)}^{(\mu^+, \mu^-)} \otimes
(F_{(m)}^{\mu^+})^{U_m} \otimes  ( F_{(\ell)}^{\mu^-})^{U_{\ell}}
$$
is a multiplicity-free sum of representations $F_{(n)}^{(\mu^+,
\mu^-)}$ of $GL_n$. Again, this can be embedded as a
$(0,1)$-subalgebra of the coordinate ring of $GL_n / U_n$.

Now repeat this with $m'$ in place of $m$ and $\ell'$ in place of
$\ell$. We again get a multiplicity-free sum of a family of
representations of $GL_n$. If we take the tensor product of the
two sums, and look at highest weight vectors for $GL_n$, we will
get a $(0,1)$-subalgebra of the tensor algebra for $GL_n$:
$$
\left ( ( \calp (M_{n, m} \oplus M_{\ell, n})/I(\widetilde{\mathcal J}_{n,m, \ell}^+ ) )^{U_m \times U_{\ell}} \otimes
( \calp (M_{n, m'} \oplus M_{\ell', n})/I(\widetilde{\mathcal J}_{n,m',\ell'}^+))^{U_{m'} \times U_{\ell'}}
\right ) ^{U_n}
$$
$$
\simeq  \bigoplus_{\mu^+, \mu^-, \nu^+, \nu^-} ( F_{(n)}^{(\mu^+, \mu^-)} \otimes  F_{(n)} ^{(\nu^+, \nu^-)} )^{U_n} \otimes
\left ( ( F_{(m)}^{\mu^+} )^{U_m} \otimes  ( F_{(\ell)}^{\mu^-})^{U_{\ell}} \otimes
(F_{(m')}^{\nu^+})^{U_{m'}} \otimes  ( F_{(\ell')}^{\nu^-})^{U_{\ell'}}.
 \right )
$$
where the sum is over partitions $\mu^+$, $\mu^-$, $\nu^+$ and
$\nu^-$ such that $\ell (\mu^+) \leq m$, $\ell (\mu^-) \leq \ell
$, $\ell (\mu^+)+l(\mu^-) \leq n$, $\ell (\nu^+) \leq m'$, $\ell
(\nu^-) \leq \ell'$, $\ell (\nu^+)+ \ell (\nu^-) \leq n$.

\medskip

On the other hand,

$$
\begin{aligned}
&\left ( ( \calp (M_{n, m} \oplus M_{\ell, n})/I(\widetilde{\mathcal J}_{n,m, \ell}^+ ) )^{U_m \times U_{\ell}} \otimes
( \calp (M_{n, m'} \oplus M_{\ell', n})/I(\widetilde{\mathcal J}_{n,m',\ell'}^+))^{U_{m'} \times U_{\ell'}}
\right ) ^{U_n} & \\
&\simeq \left ( (\calp (M_{n, m} \oplus M_{\ell,
n})/I(\widetilde{\mathcal J}_{n,m, \ell}^+ )) \otimes  (
\calp (M_{n, m'} \oplus M_{\ell', n})/I(\widetilde{\mathcal
J}_{n,m',\ell'}^+))\right)^{U_n \times U_m \times U_{\ell} \times U_{m'} \times U_{\ell'}} & \\
&\simeq \left ( \calp (M_{n, m} \oplus M_{\ell,n} \oplus M_{n, m'} \oplus M_{\ell', n})/I(\widetilde{\mathcal{J}}^+_{n,m,\ell} \oplus \widetilde
{\mathcal{J}}^+_{n,m',\ell'}) \right)^{U_n \times U_m \times U_{\ell} \times U_{m'} \times U_{\ell'}} & \\
&\simeq \left( \calp (M_{n, m+m'} \oplus M_{\ell +\ell',n}) /I(\widetilde{\mathcal{J}}^+_{n,m,\ell} \oplus \widetilde
{\mathcal{J}}^+_{n,m',\ell'}) \right)^{U_n \times U_m \times U_{\ell} \times U_{m'} \times U_{\ell'}} & \\
&\simeq \left( \bigoplus_{\la^+, \la^-} (F_{(n)}^{(\la^+,
\la^-)})^{U_n} \otimes \widetilde F_{(m+m',\ell+\ell')}^{(\la^+,\la^-)}
/ I(\widetilde{\mathcal{J}}^+_{n,m,\ell} \oplus \widetilde
{\mathcal{J}}^+_{n,m',\ell'}) \right)^{U_m \times U_{\ell} \times U_{m'} \times U_{\ell'}} & \\
&\simeq \bigoplus_{\la^+, \la^-} (F_{(n)}^{(\la^+,\la^-)})^{U_n}
\otimes \left( \widetilde F_{(m+m',\ell+\ell')}^{(\la^+,\la^-)}
/ (\widetilde{\mathcal{J}}^+_{n,m,\ell} \oplus \widetilde
{\mathcal{J}}^+_{n,m',\ell'}) \cdot \widetilde F_{(m+m',\ell+\ell')}^{(\la^+,\la^-)} \right)^{U_m \times U_{\ell} \times U_{m'} \times U_{\ell'}} &
\end{aligned}
$$
which tells us about the $\gl_{m ,\ell} \oplus \gl_{m',\ell'}$ decomposition in the representation
$\widetilde F_{(m+m', \ell+\ell')}^{(\lambda^+, \lambda^-)}$ of $\gl_{(m+m'), (\ell + \ell')}$.
This completes the proof. $\qquad \square$

\vskip 10pt

The construction of \S 4.1 of course is just the case $\ell = 0 =
\ell'$ of the current discussion. That case is notable for staying
completely in the context of finite-dimensional representation
theory.  Another case of interest is when $\ell = 0  = m'$. Then,
although the representations of $\gl_{m , \ell'}$ are infinite
dimensional, the representations of the subalgebras $\gl_m$ and
$\gl_{\ell'}$ are finite dimensional. This case is analogous to
branching from $GL_n$ to $O_n$ (or from $GL_{2n}$ to $Sp_{2n}$).

\section{The Stable Range and Relations Between Reciprocity Algebras}

Let us summarize our discussions this far. Given any classical
symmetric pair, we can embed it in a (family of) see-saw symmetric
pair(s). Doing this, we find that (a $(0,1)$-subalgebra of) the
branching algebra for the pair can equally well be interpreted as
the branching algebra for a dual family of representations of the
dual symmetric pair. The representations of the dual symmetric
pair will frequently be infinite dimensional, but they are always
highest weight modules.

An immediate consequence of this isomorphism of algebras is the isomorphisms of intertwining spaces and hence equality of multiplicities, which we have collectively described as {\it reciprocity laws}. These reciprocity laws are of the same nature as Frobenius Reciprocity for induced representations of groups.

From \S 4.2, we see that the see-saw symmetric pairs actually come
in two parameter families. If one of the pairs involves many more
variables than the other, then certain features of the discussions
above become simpler.

Take the results of Theorem 4.1 as an illustration: The
Littlewood-Richardson coefficients for $GL_n$,
$$
c_{\mu \nu}^\lambda = \dim \Hom_{GL_n} ( F_{(n)}^\lambda , F_{(n)}^\mu \otimes  F_{(n)}^\nu)
$$
$$
=  \dim \Hom_{GL_m \times GL_{\ell} } ( F_{(m)}^\mu \otimes
F_{(\ell)}^\nu, F_{(m + \ell)}^\lambda)
$$
are independent of $n$, if $n \geq m+ \ell$, and depend only on the
shape of the partitions $\mu$, $\nu$ and $\lambda$.

Consider another example: branching from $GL_{2n}$ to $Sp_{2n}$. If we let these groups act on $\calp (M_{2n,m})$, we get the see-saw pairs $(Sp_{2n},\so_{2m})$ and $(GL_{2n}, \gl_m)$. The branching coefficients $d_{\lambda}^\mu$ from $GL_{2n}$ to $Sp_{2n}$ can be described as follows:
$$
F_{(2n)}^\lambda \mid_{Sp_{2n}} = \sum_{\mu} d_{\lambda}^\mu \ V_{(2n)}^\mu
$$
where
$$
\begin{aligned}
d_{\lambda}^\mu &= \dim \Hom_{Sp_{2n}}(V^\mu_{2n}, F^\la_{2n}) \\
&= \dim \Hom_{GL_m} \left( F_{(m)}^\lambda, F_{(m)}^\mu
\otimes  \calS ( \so_{2m}^{(2,0)}) \right) \\
&= \dim \Hom_{GL_m} \left( F_{(m)}^\lambda, F_{(m)}^\mu
\otimes  \calS ( \wedge^2 \C^m) \right)
\end{aligned}
$$
is independent of $n$, if $n \geq m$, and only depends on the
diagrams $\lambda$ and $\mu$.  This allows one to create a theory
of ``stable characters'' for $Sp_{2n}$. Similar considerations
apply to $GL_n$ and $O_n$ and this idea has been actively pursued by
\cite{koike-terada}, amongst others.

These are all instances of stability laws. The well-known one-step branching from $GL_n$ to $GL_{n-1}$ is another instance. This branching can be described entirely by diagrams, with no mention of the size $n$, if $n$ is large. Iteration of this branching also shows that when $n$ is large, the weight multiplicities of dominant weights of an irreducible $GL_n$ representation are independent of $n$. See \cite{benkart} for the other classical groups, which don't share this stability property.

In the sections that follow, we will illustrate the simplifications
that occur in the stable range, highlighting certain specific
see-saw pairs. In all these cases, we show that the branching
algebras associated to symmetric pairs can all be described by use
of suitable branching algebras associated to the general linear
groups.  Thus, if we can have control of the solution in the general
linear group case, we will have some control of the other classical
groups.  The other non-trivial examples will be important extensions
of this work, and we hope to see them in further papers, for
example, \cite{HTW1}, \cite{HTW2} and \cite{howe-lee}.

\section{Stability for Branching from $GL_n$ to $O_n$}

We begin with a detailed discussion of the case of $(GL_n, O_n)$
and $(\spm_{2m}, GL_m)$. Here we have already encountered the
stable range, without the name. It is when $n> 2m$. Several things
happen in the stable range:
\begin{enumerate}

\item[(a)] The representations $E_{(n)}^\mu$ of the orthogonal
group remain irreducible when restricted to the special orthogonal
group $SO_n$, and furthermore, no two of them are equivalent.

\item[(b)] Recall the algebra ${\mathcal J}_{n,m}$  of
$O_n$-invariant polynomials on $M_{n, m}$ generated by the
quadratic invariants, which is the abelian subalgebra $\spm_{2m}^{(2, 0)}$ of
$\spm_{2m}$. In the stable range (in fact it holds true whenever $n \geq m$), the natural surjective homomorphism
$$
\calS (\spm_{2m}^{(2, 0)}) \rightarrow {\mathcal J}_{n,m}
$$
is an isomorphism. See Theorem 3.3(a).

\item[(c)] In the stable range, the multiplication map
$$
\Har_{n,m} \otimes  {\mathcal J}_{n,m} \simeq \Har_{n,m} \otimes  \calS (\spm_{2m}^{(2, 0)}) \rightarrow \calp (M_{n,m})
$$
is also an isomorphism of $O_n \times GL_m$-modules.  See Theorem
3.3(a).

\end{enumerate}

Of course, the subspace $\Har_{n,m}$ of harmonic polynomials is not an algebra -- it
is not closed under multiplication. This is quite clear, since $\Har_{n,m}$ contains
all the linear functions, which generate the whole polynomial ring.
However, to form the reciprocity algebra associated to the symmetric see-saw
pairs $(GL_n, O_n)$ and $(\spm_{2m}, GL_m)$, we need to take the $U_{O_n}$-invariants.
Thus, our reciprocity algebra is a subalgebra of
$$
(\Har_{n,m} \otimes   \calS (\spm_{2m}^{(2, 0)}) )^{U_{O_n}} =
\Har_{n,m}^{U_{O_n}} \otimes   \calS (\spm_{2m}^{(2, 0)})
$$
$$
\simeq \left (\bigoplus_\mu (E_{(n)}^\mu )^{U_{O_n}} \otimes  F_{(m)}^\mu \right ) \otimes
 \calS (\spm_{2m}^{(2, 0)}). \eqno (9.1)
$$

\begin{thm} When $n > 2m$, the space  $\Har_{n,m}^{U_{O_n}}$ is
a subalgebra of $\calp (M_{n,m})$. Hence, the algebra
$\calp (M_{n, m})^{U_{O_n}}$ is isomorphic to a tensor product
$$
\calp (M_{n, m})^{U_{O_n}} \simeq \Har_{n,m}^{U_{O_n}} \otimes   \calS (\spm_{2m}^{(2, 0)})
$$
of the algebras $\Har_{n,m}^{U_{O_n}}$ and $\calS (\spm_{2m}^{(2, 0)})$. Furthermore, the
algebra  $\Har_{n,m}^{U_{O_n}}$ is isomorphic (as a representation) to the subalgebra
$\calr^+(GL_m/U_m)$ of $\calr (GL_m/U_m)$ defined
by the polynomial representations.
\end{thm}

\nt {\bf Proof.} Note that $\Har_{n,m}^{U_{O_n}}$ can be identified
with a subalgebra $\calr^+(GL_m/U_m)$ of $\calr (GL_m/U_m)$ defined
by the polynomial representations, from our discussion in \S 3.2.
Consider the space of polynomials belonging to the sum in the last
expression of equation (9.1). Let \\
$\{ x_{jk} \mid j=1,\ldots,n, k=1,\ldots,m \}$ be the standard
matrix entries on $M_{n, m}$. In order to make the unipotent group
$U_{O_n}$ of $O_n$ maximally compatible with (in fact, contained in)
the unipotent subgroup $U_n$ of $GL_n$, we choose the inner product
on ${\C}^n$ as in Section 3.3.1.  By this choice, joint $O_n \times
GL_m$ harmonic highest weight vectors are monomials in the
determinants
$$
\delta_j = \det \left[ \begin{matrix} x_{11} &x_{12} &\ldots
&x_{1j} \cr x_{21} &x_{22} &\ldots &x_{1j} \cr \vdots &\vdots
&\vdots &\vdots \cr x_{j1} &x_{j2} &\ldots &x_{jj} \end{matrix}
\right] \quad \text{ for }  j=1,\ldots, m.\quad \eqno (9.2)
$$
From this, we can see that the space
$\sum_\mu ( E_{(n)}^\mu )^{U_{O_n}} \otimes  F_{(m)}^\mu$ is spanned by the
monomials in the determinants
$$
\det \left[ \begin{matrix} x_{1,b_1} &x_{1,b_2} &\ldots
&x_{1,b_j}\cr x_{2,b_1} &x_{2,b_2} &\ldots &x_{2,b_j} \cr \vdots
&\vdots    &\vdots &\vdots \cr x_{j,b_1}  &x_{j,b_2} &\ldots
&x_{j,b_j} \end{matrix} \right] \eqno (9.3)
$$
as $\{ b_1, b_2, b_3, \ldots , b_j \}$ ranges over
all $j$-tuples of integers from 1 to $m$. Indeed, the span of such monomials
is clearly invariant under $\gl_m$, and consists of highest weight vectors for
$O_n$. Finally, we see that these monomials will all be harmonic, because the
partial Laplacians spanning $\spm_{2m}^{(0,2)}$ have the form
$$
\Delta_{ab} = \sum_{j = 1}^n \frac{\partial^2}{\partial x_{j,a}
\partial x_{n+1-j, b}}.\eqno (9.4)
$$
Since every term of $\Delta_{ab}$ involves differentiating with
respect to a variable $x_{jk}$ with $j > n/2$, and the
determinants (9.3) do not depend on these variables, we see that
they will be annihilated by the $\Delta_{ab}$, which means that
they are harmonic.  This shows that $\mcH^{U_{O_n}}_{n,m}$ is a
subalgebra of $\mcP(M_{n,m})$.

We have thus completed the proof of the theorem. $\quad \square$

\vskip 10pt

We can use the description in Theorem 9.1 of $\calp (M_{n,
m})^{U_{O_n}}$ to relate the branching algebra $ \calp (M_{n,
m})^{U_{O_n} \times U_m}$ to the tensor product algebra for
$GL_m$. As a $GL_m$-module, the space $\spm_{2m}^{(2, 0)}$ is
isomorphic to $\calS^2 ({\C}^m)$, the space of symmetric $m \times
m$ matrices. It is well known that the symmetric algebra $\calS (
\calS^2 ({\C}^m))$ is multiplicity-free as a representation of
$GL_m$, and decomposes into a sum of one copy of each polynomial
representation corresponding to a diagram with rows of even length
(or a partition of even parts):
$$
\calS (\calS^2 ({\C}^m)) \simeq \bigoplus_\nu F_{(m)}^{2\nu}.
\eqno (9.5)
$$
(Note that this result is in several places in the literature. See
\cite{goodman-wallach-book} and \cite{howe-schur} for example.)

As a $GL_m$-module, $\calS (\calS^2 ({\C}^m))$ could be embedded in
$\calr (GL_m/U_m)$, but the algebra structures on these two algebras
are quite different.

Using the dominance filtration (see \S 3.2), we have a canonical
$\widehat A^+$-algebra filtration on $\mcS(\mcS^2({\bbC}^m))$. If we
form the associated graded algebra, then Theorem 3.2  says that it
will be isomorphic to the subalgebra of $R(GL_m /U_m)$ spanned by
the representations attached to diagrams with even length rows.

Let us denote the associated graded algebra of $\calS (\calS^2
({\C}^m))$ by $\text{Gr}_{\widehat A_m^+} \calS (\calS^2 ({\C}^m))$.
Let us denote the subalgebra of $\calr (GL_m/U_m)$ spanned by the
representations attached to diagrams with even length rows by
$\calr^{+2}(GL_m/U_m)$.

We can filter the tensor product $\Har_{n,m}^{U_{O_n}} \otimes
\calS (\spm_{2m}^{(2, 0)})$ by means of the filtration on  $ \calS
(\spm_{2m}^{(2, 0)})$. The associated graded algebra will then be
$\Har_{n,m}^{U_{O_n}} \otimes\, \text{Gr}_{\widehat A_m^+} \calS
(\spm_{2m}^{(2, 0)})$. This discussion has indicated that the
following result holds.

\vskip 10pt

\begin{thm} When $n > 2m$, the associated graded algebra of
$\calp (M_{n,m})^{U_{O_n}}$ with respect to the dominance filtration on the factor
${\mathcal J}_{n,m}$  is isomorphic to the tensor product of the graded subalgebras
$\calr^+(GL_m/U_m)$ and $\calr^{+2}(GL_m/U_m)$ of $\calr (GL_m /U_m)$:
$$
\text{Gr}_{\widehat A_m^+} (\calp (M_{n,m})^{U_{O_n}}) \simeq \calr^+(GL_m/U_m) \otimes  \calr^{+2}(GL_m/U_m).
$$
\end{thm}

\vskip 10pt

Of course, $\text{Gr}_{\widehat A_m^+} (\calp (M_{n,m})^{U_{O_n}})$
is isomorphic as a $GL_m$-module to $ \calp (M_{n,m})^{U_{O_n}}$ in
an obvious way, by construction.  Also $\text{Gr}_{\widehat A_m^+}
(\calp (M_{n,m})^{U_{O_n}})$ inherits the $\widehat A_{O_n}^+$
grading from $ \calp (M_{n,m})^{U_{O_n}}$ -- it becomes identified
with the $\widehat A_m^+$ grading on the first factor $\calr^+(GL_m
/U_m)$ in the tensor product of Theorem 9.2. On the other hand, the
second factor is also $\widehat A_m^+$-graded in the obvious way,
since it is the factor which defines the associated graded. When we
take the $U_m$ invariants, we get another grading by $\widehat
A_m^+$, associated to the $A_m$ action on the $U_m$ invariants. This
triply $\widehat A_m^+$-graded algebra is evidently a
$(0,1)$-subalgebra of the tensor product algebra of $GL_m$.

On the other hand, we could take the $U_m$ invariants inside
$\calp (M_{n,m})^{U_{O_n}}$, and then pass to the associated graded. It is not
hard to convince oneself that these two processes commute with each other.
Hence, we finally have:

\vskip 10pt

\nt {\bf Corollary 9.3} {\it When $n > 2m$, the associated graded
algebra of $U_m$ invariants in $\calp (M_{n,m})^{U_{O_n}}$,
$$
\text{Gr}_{\widehat A_m^+} \left( \left( \calp (M_{n,m})^{U_{O_n}} \right)^{U_m} \right)
\simeq \left( \text{Gr}_{\widehat A_m^+} ( \calp (M_{n,m})^{U_{O_n}} ) \right)^{U_m}
$$
$$
\simeq \left( \calr^+(GL_m/U_m) \otimes  \calr^{+2}(GL_m/U_m) \right)^{U_m}
$$
is a triply-graded $(0,1)$-subalgebra of the tensor product
algebra of $GL_m$. The restrictions on the gradings which define
$\text{Gr}_{\widehat A_m^+} \left( \left( \calp
(M_{n,m})^{U_{O_n}} \right)^{U_m} \right)$ are:
\begin{enumerate}
\item[(a)] the weight on the first factor of
$\left(\calr (GL_m/U_m) \otimes  \calr (GL_m/U_m) \right)^{U_m}$ should correspond to
a partition (i.e., it should be a polynomial weight), and
\item[(b)] the weight on the
second factor should correspond to a partition with even parts.
\end{enumerate}
}

\vskip 10pt

\nt {\bf Remark:} The content of Corollary 9.3 in terms of
multiplicities is the Littlewood Restriction Formula
\cite{ew-annals}, \cite{HTW1}; see formula (2.4.1),
\cite{king-modify}; see (5.7) with (4.19), \cite{koike-terada};
see Theorem 1.5.3 and 2.3.1, \cite{littlewood-paper} and
\cite{littlewood}. With this result it is possible to compute a
basis of the reciprocity algebra for $(GL_n, O_n)$ using
\cite{HTW1}; see second preprint of \cite{howe-lee}.

\bigskip

\section{Tensor Products for $O_n$}

According to Theorem 6.1, we can compute tensor products for the
orthogonal group via the algebra
$$
\left( \left( \calp (M_{n,m})/I({\mathcal J}_{n,m}^+) \right )^{U_m} \otimes
\left ( \calp (M_{n,\ell})/I({\mathcal J}_{n,\ell}^+) \right )^{U_{\ell}} \right ) ^{U_{O_n}}.
$$
Here the stable range is $n > 2(m + \ell) $. Then we have
$$
\calp (M_{n, m + \ell}) \simeq \Har_{n,m+\ell} \otimes  \calS
(\spm_{2(m + \ell)} ^{(2,0)}).
$$
Furthermore,
$$
\calS (\spm_{2(m + \ell)} ^{(2,0)}) = \calS (\spm_{2m} ^{(2,0)} \oplus
\spm_{2 \ell} ^{(2,0)} \oplus ({\C}^{m} \otimes  {\C}^{\ell}) )
$$
$$
\qquad \qquad \qquad \simeq \calS (\spm_{2m} ^{(2,0)}) \otimes
\calS (\spm_{2\ell} ^{(2,0)}) \otimes \calS ({\C}^{m} \otimes
{\C}^{ \ell})
$$
Since ${\mathcal J}_{n,m} \simeq \calS (\spm_{2m} ^{(2,0)})$ and
${\mathcal J}_{n,\ell} \simeq \calS (\spm_{2 \ell}^{(2,0)}) $, we see that
$$
\calp (M_{n,m})/I({\mathcal J}_{n,m}^+)  \otimes
 \calp (M_{n,\ell})/I({\mathcal J}_{n,\ell}^+) \simeq \calp (M_{n,m} \oplus M_{n,\ell})/I({\mathcal J}_{n,m}^+
 \oplus {\mathcal J}_{n,\ell}^+)
$$
$$
\qquad \qquad \qquad \simeq \Har_{n,m+\ell} \otimes \calS
({\C}^{m} \otimes  {\C}^{\ell}). \eqno (10.1)
$$
Thus, using equation (10.1), we see that
$$
\left( \calp (M_{n,m})/I({\mathcal J}_{n,m}^+)  \otimes
 \calp (M_{n,\ell})/I({\mathcal J}_{n,\ell}^+) \right) ^{U_{O_n}}
$$
$$
\simeq \left( \Har_{n,m + \ell}  \otimes  \calS ({\C}^{m} \otimes  {\C}^{ \ell}) \right) ^{U_{O_n}}
\simeq \Har_{n,m+\ell}^{U_{O_n}} \otimes  \ \calS ({\C}^{m} \otimes  {\C}^{ \ell})
$$
$$
\simeq \left ( \bigoplus_\lambda E_{(n)}^\lambda \otimes  F_{(m + \ell)}^\lambda \right )^{U_{O_n}} \otimes  \ \calS ({\C}^{m} \otimes  {\C}^{ \ell})
$$
$$
\simeq  \left( \bigoplus_\lambda ( E_{(n)}^\lambda )^{U_{O_n}}
\otimes  F_{(m + \ell)}^\lambda \right )
\otimes  \calS ({\C}^{m} \otimes  {\C}^{ \ell})
$$
$$
\simeq  \left( \bigoplus_\lambda ( F_{(n)}^\lambda )^{U_n}
\otimes  F_{(m + \ell)}^\lambda \right )
\otimes  \calS ({\C}^{m} \otimes  {\C}^{ \ell}).
$$
Note that $F_{(n)}^\lambda$ is the $GL_n$ representation generated by the highest weight of the $O_n$ representation $E_{(n)}^\lambda$ and both $(F_{(n)}^\lambda)^{U_n }$ and $(E_{(n)}^\lambda)^{O_n}$ are one dimensional.

Hence, finally we get
$$
\left ( \left ( \calp (M_{n,m})/I({\mathcal J}_{n,m}^+)  \otimes
 \calp (M_{n,\ell})/I({\mathcal J}_{n,\ell}^+) \right ) ^{U_{O_n}} \right ) ^{U_m \times U_{\ell}}
$$
$$
\simeq \left (\left( \bigoplus_\lambda ( F_{(n)}^\lambda )^{U_n}
\otimes  F_{(m + \ell)}^\lambda  \right ) \otimes  \calS ({\C}^{m}
\otimes  {\C}^{ \ell}) \right )^{U_m \times U_{\ell}}. \eqno
(10.2)
$$

We can interpret this algebra in term of tensor product algebras
for general linear groups. According to Theorem 3.4(c) and (3.8),
as a $GL_m \times GL_{\ell}$ module, we have
$$
\calS ({\C}^{m} \otimes  {\C}^{ \ell})  \simeq \bigoplus_\delta
F_{(m)}^\delta  \otimes  F_{(\ell)}^\delta. \eqno (10.3)
$$
We also know that
$$
\calr (GL_m /U_m \times GL_{\ell} / U_{\ell}) \simeq
\bigoplus_{\mu, \nu} F_{(m)}^\mu \otimes  F_{(\ell)}^\nu. \eqno
(10.4)
$$
Since this algebra is bigraded by $\mu$ and $\nu$, we can consider the \lq \lq diagonal" $(0,1)$-subalgebra
$$
\Delta \calr (GL_m /U_m \times GL_{\ell} / U_{\ell}) \simeq
\bigoplus_{\delta} F_{(m)}^\delta \otimes  F_{(\ell)}^\delta \eqno
(10.5)
$$
resulting from requiring the two partitions to be the same.
Evidently, the algebra $\calS ({\C}^m \otimes  {\C}^{\ell})$ is
isomorphic to $\Delta \calr (GL_m /U_m \times GL_{\ell} / U_{\ell})$
as $GL_m \times GL_{\ell}$-module. They are not isomorphic as
algebras, since $\Delta \calr (GL_m /U_m \times GL_{\ell} /
U_{\ell})$ is graded by $\widehat A_m^+ \times \widehat A_\ell^+$,
while $\calS ({\C}^m \otimes  {\C}^{\ell})$ is not. However, we may
filter $\calS ({\C}^m \otimes  {\C}^{\ell})$ by the representations
of $GL_m \times GL_{\ell}$ (or of either factor) using the dominance
filtration (see \S 3.2), and then the associated graded algebra will
be isomorphic to $\Delta \calr (GL_m /U_m \times GL_{\ell} /
U_{\ell})$ by Theorem 3.2:
$$
\text{Gr}_{\widehat A_m^+ \times \widehat A_\ell^+} (\calS ({\C}^m
\otimes  {\C}^{\ell})) \simeq \Delta \calr (GL_m /U_m \times
GL_{\ell} / U_{\ell}).
$$

Now turn to the first factor $\bigoplus_\lambda \left( (
F_{(n)}^\lambda )^{U_n}  \otimes  F_{(m + \ell)}^\lambda \right)$
on the right hand side of equation (10.2).  According to Theorem
4.1, we can write this as
$$
\bigoplus_\lambda ( F_{(n)}^\lambda )^{U_n}  \otimes  F_{(m +
\ell)}^\lambda  \simeq \bigoplus_{\alpha, \beta } \left (
F_{(n)}^\alpha \otimes  F_{(n)}^\beta \right )^{U_n} \otimes
F_{(m)}^\alpha \otimes  F_{(\ell)}^\beta . \eqno (10.6)
$$
Combining equations (10.2), (10.3) and (10.6), we see that
$$
\left ( \left ( \calp (M_{n,m})/I({\mathcal J}_{n,m}^+)  \otimes
 \calp (M_{n,\ell})/I({\mathcal J}_{n,\ell}^+)
 \right )^{U_{O_n}} \right ) ^{U_m \times U_{\ell} }
$$
$$
\qquad \simeq
\left ( \left ( \bigoplus_{\alpha, \beta} \left ( F_{(n)}^\alpha \otimes
F_{(n)}^\beta \right )^{U_n} \otimes  F_{(m)}^\alpha \otimes  F_{(\ell)}^\beta \right ) \otimes
\left ( \bigoplus_{\delta } F_{(m)}^\delta
\otimes  F_{(\ell)}^\delta \right ) \right )^{U_m \times U_{\ell}}\qquad
$$
$$
\qquad \simeq  \left ( \bigoplus_{\alpha, \beta, \delta} \left ( F_{(n)}^\alpha \otimes  F_{(n)}^\beta \right
)^{U_n} \otimes  \left ( F_{(m)}^\alpha \otimes  F_{(m)}^\delta \right ) \otimes
\left ( F_{(\ell)}^\beta \otimes  F_{(\ell)}^\delta \right )
  \right ) ^{U_m \times U_{\ell}} \eqno (10.7)
$$
$$
\simeq \bigoplus_{\alpha, \beta, \delta } \left( F_{(n)}^\alpha
\otimes  F_{(n)}^\beta  \right)^{U_n} \otimes \left(
F_{(m)}^\alpha \otimes  F_{(m)}^\delta \right)^{U_m} \otimes
\left( F_{(\ell)}^\beta \otimes F_{(\ell)}^\delta
\right)^{U_{\ell}} \qquad\qquad
$$
At this point, this is an isomorphism of graded vector spaces,
not an algebra isomorphism.

We may interpret the last expression in (10.7) analogously to
(10.4) and (10.5). We have the (polynomial) tensor product
algebras
$$
(\calr^+(GL_k/ U_k) \otimes  \calr^+(GL_k / U_k) )^{U_k} \simeq
\bigoplus_{\lambda, \mu} \left ( F_{(k)}^\lambda \otimes  F_{(k)}^\mu \right )^{U_k}
$$
for $k = n, m$ and $\ell$. If we form the tensor product of these, we get
$$
(\calr^+(GL_n/ U_n) \otimes  \calr^+(GL_n / U_n) )^{U_n} \otimes
(\calr^+(GL_m/ U_m) \otimes  \calr^+(GL_m / U_m) )^{U_m}
$$
$$
\otimes  \ (\calr^+(GL_{\ell}/ U_{\ell}) \otimes  \calr^+(GL_{\ell} / U_{\ell}) )^{U_{\ell}}
$$
$$
\simeq \bigoplus_{\alpha, \beta, \delta, \lambda, \mu, \nu}
\left (F_{(n)}^\alpha \otimes  F_{(n)}^\beta \right )^{U_n} \otimes
\left (F_{(m)}^\delta \otimes  F_{(m)}^\lambda \right )^{U_m} \otimes
\left (F_{(\ell)}^\mu \otimes  F_{(\ell)}^\nu \right )^{U_{\ell}}
$$
Let us denote this algebra by $\TT_{n, m, \ell}$. The algebra
$\TT_{n, m, \ell}$ is $(\widehat A_n^+)^3 \times (\widehat
A_m^+)^3 \times (\widehat A_{\ell}^+)^3$-graded. If we require
that $\lambda = \alpha$, or that $\mu = \beta$, or that $\nu =
\delta$, then we obtain $(0,1)$-subalgebras of $\TT_{n, m, \ell}$.
If $\delta = \alpha$, we will denote it by $\Delta_{1,3}
\TT_{n,m,\ell}$, and so forth. The subalgebra obtained by
requiring that all three diagonal conditions occur at once will be
denoted by using all three $\Delta$'s. Thus we will write
$$
\Delta_{1,3} \Delta_{2,5} \Delta_{4,6} \TT_{n, m, \ell}
\qquad\qquad \qquad \qquad \qquad \qquad \qquad
$$
$$
=\sum_{\alpha, \beta, \delta} \left ( F_{(n)}^\alpha \otimes
F_{(n)}^\beta  \right )^{U_n} \otimes  \left ( F_{(m)}^\alpha
\otimes  F_{(m)}^\delta \right )^{U_m} \otimes \left (
F_{(\ell)}^\beta \otimes  F_{(\ell)}^\delta \right ) ^{U_{\ell}}
\eqno (10.8)
$$

We see from equations (10.7) and (10.8), that $\Delta_{1,3}
\Delta_{2,5} \Delta_{4,6}  \TT_{n, m, \ell}$ and
\[ \left ( \left ( \calp (M_{n,m})/I({\mathcal J}_{n,m}^+)  \otimes
 \calp (M_{n,\ell})/I({\mathcal J}_{n,\ell}^+) \right ) ^{U_{O_n}} \right ) ^{U_m \times
U_{\ell} }\]
are isomorphic as multigraded vector spaces. They may
not be isomorphic as algebras, because
$$
\left ( \left ( \calp (M_{n,m})/I({\mathcal J}_{n,m}^+)  \otimes
 \calp (M_{n,\ell})/I({\mathcal J}_{n,\ell}^+) \right ) ^{U_{O_n}} \right ) ^{U_m \times
U_{\ell} }
$$
is not graded, while we see from equations (10.7) and (10.8), that
$\Delta_{1,3} \Delta_{2,5} \Delta_{4,6} \TT_{n, m, \ell}$ is.
However, if we pass to the associated graded of $\calS ({\C}^m
\otimes  {\C}^{\ell})$, then the two algebras do become isomorphic.
We record this fact.

\vskip 10pt

\begin{thm} Assume the stable range $n > 2(m+\ell)$.  We have the following isomorphisms of $(\widehat A_n^+)^3 \times (\widehat A_m^+)^3 \times (\widehat A_{\ell}^+)^3$-graded algebras:
$$
\text{Gr}_{(\widehat A_n^+)^3 \times (\widehat A_m^+)^3 \times (\widehat A_{\ell}^+)^3} \left( \left( \left( \calp (M_{n,m})/I({\mathcal J}_{n,m}^+ )  \otimes
 \calp (M_{n,\ell})/I({\mathcal J}_{n,\ell}^+) \right) ^{U_{O_n}} \right)^{U_m \times
U_{\ell} } \right)
$$
$$
\qquad \qquad \qquad \qquad \simeq \Delta_{1,3} \Delta_{2,5}
\Delta_{4,6} \TT_{n, m, \ell}.
$$
\end{thm}

\nt {\bf Remark:} The content of Theorem 10.1 in terms of
multiplicities can be found in \cite{HTW1}; see formula (2.1.2),
\cite{king-s-fcn}; see Theorem 4.1 and \cite{newell}.

\section{Restriction from $O_{n+m}$ to $O_n \times O_{m}$}

Consider now branching from $O_{n+m}$ to the product $O_n \times
O_{m}$. We look at the action of $O_{n+m}$ on $M_{n+m, \ell}$ by
multiplication on the left.  Here we simply state the relationship
of the reciprocity algebra for $(O_{n+m},O_n \times O_m)$ in the
stable range, in relation to the $GL_n$ tensor product algebras.
We shall omit the proof.

The content of Theorem 11.1 in terms of multiplicities can be
found in \cite{HTW1}; see formula (2.2.2), \cite{king-s-fcn}; see
(2.16) and \cite{koike-terada}; see Theorem 2.5 and corollary 2.6.
With this result it is possible to compute a basis of the
reciprocity algebra for $(O_{n+m}, O_n \times O_m)$ using
\cite{HTW2}; see second preprint of \cite{howe-lee}.

\begin{thm} Assume the stable range $\min\, (m,n) > 2\ell$.  We have the following isomorphisms of $(\widehat A_\ell^+)^3 \times (\widehat A_\ell^+)^3 \times (\widehat A_{\ell}^+)^3$-graded algebras:
$$
\text{Gr}_{(\widehat A_\ell^+)^3 \times (\widehat A_\ell^+)^3
\times (\widehat A_{\ell}^+)^3} \left( \calp (M_{n+m ,
\ell})/I({\mathcal J}_{n+m,\ell}^+) \right) ^{U_{O_n}\times
U_{O_m} \times U_{\ell} }
$$
$$
\qquad \qquad \simeq
\left( \calr^+(GL_{\ell}/U_{\ell}) \otimes  \calr^+(GL_{\ell}/U_{\ell}) \otimes  \calr^{+2}(GL_{\ell} /U_{\ell}) \right)^{U_{\ell}}.
$$
\end{thm}

\newpage

\nt {\bf Appendix: A Proof of the Separation of Variables Theorem}

\vskip 10pt

We provide a simple proof here for $G = O_n$. It could be adapted easily for $Sp_{2n}$ acting on copies of $\C^{2n}$ or $GL_n$ acting on copies of $\C^n$ and $\C^{n*}$.

Let $O_n$ act by the usual left multiplication on $M_{n,m}$, the
$n \times m$ matrices. We assume that $n \geq 2m$, which includes
the stable range $n > 2m$. Let $r_{ij}$ be the invariant pairing
between the $i$-th column and the $j$-th column. Recall that
$\mathcal J_{n,m} = \calp (M_{n,m})^{O_n}$ is generated freely by
the homogeneous quadratic polynomials $\{ r_{ij} \}$ in this range
(see \cite{goodman-wallach-book} Theorem 5.2.7). Also, the space
of harmonics, i.e., polynomials annihilated by all the
differential operators dual to the $r_{ij}$'s, is denoted by
$\calh_{n,m}$ as in Theorem 3.3(a).

\vskip 5pt

\nt {\bf Theorem (Separation of Variables)} {\it If $n \geq 2m$,
then}
$$
\calp (M_{n,m}) \simeq \calh_{n,m} \otimes \calp (M_{n,m})^{O_n}
$$
\vskip 5pt

\nt {\bf Remarks.}  Proofs of this result for orthogonal groups (see
Theorem 2.5 of \cite{ton-that-76}) and for symplectic groups (see
Theorem 1.10 of \cite{ton-that-77}) are given by Ton-That using
results of \cite{kostant}.

\vskip 5pt

\nt {\bf Proof of Separation of Variables Theorem:}  Let $I (
{\mathcal J_{n,m}^+})$ be the ideal in $\calp (M_{n,m})$ generated
by $r_{ij}$'s with zero constant terms, and consider $\calp
(M_{n,m})$ as an $I ( {\mathcal J_{n,m}^+})$ module by
multiplication.

We shall need the notion of a regular sequence. First, denote the
ideal in $\calp (M_{n,m})$ generated by $\{ f_1, \ldots, f_s \}$
by the symbol $<f_1,\ldots,f_s>\calp (M_{n,m})$.  A sequence
$\{f_1, \ldots, f_k\} \subset  I ( {\mathcal J_{n,m}^+})$ forms a
{\it regular sequence} for $\calp (M_{n,m})$ if
\begin{enumerate}
\item[(a)] $f_i$ is not a zero-divisor on $\calp (M_{n,m})/<f_1, \ldots, f_{i-1}> \calp (M_{n,m})$ for all $i=1, \ldots, k$, and
\item[(b)] $\calp (M_{n,m})/<f_1, \ldots, f_k> \calp (M_{n,m})$ is non-zero.
\end{enumerate}

Geometrically, saying that a given function $f$ is not a
zero-divisor on \\
$\calp (M_{n,m})/<f_1, \ldots, f_{i-1}>$ is the
same as saying that $f$ does not vanish identically on any
irreducible component of the zero set of $\{ f_1, \ldots, f_{i-1}
\}$.  This in turn is the same as saying that each irreducible
component of $\{ f_1, \ldots, f_{i-1}, f \}$ has dimension one less
that the component of the zero set of $\{ f_1, \ldots, f_{i-1} \}$.

Separation of variables would follow from knowing that the $r_{ij}$'s (in some order) form a regular sequence for $\calp (M_{n,m})$. In fact, you can take any order you want. For an ideal $I$ in a commutative ring $S$, we have the chain:
$$
S \supseteq I \supseteq I^2 \supseteq I^3 \supseteq \ldots
$$
and we can thus form the associated graded algebra
$$
\text{Gr}_I S = S/I \oplus I/I^2 \oplus I^2/I^3 \oplus \ldots
$$
with multiplication (setting $S =I^0$)
$$
I^{i}/I^{i+1} \otimes I^j/I^{j+1} \longrightarrow I^{i+j}/I^{i+j+1}
$$
induced by multiplication on $S$.  In our context, $S=\calp
(M_{n,m})$ and $\calp (M_{n,m})/{I ( {\mathcal J_{n,m}^+})}$ is
the coordinate ring of the null-cone and as a linear space (in
particular, as a $O_n \times GL_m$ module), it is isomorphic to
$\calh_{n,m}$.  If $\{ r_{ij}\} \subset I ( {\mathcal J_{n,m}^+})$
is a regular sequence of $\calp (M_{n,m})$, then we have a nice
presentation of $\text{Gr}_{I ( {\mathcal J_{n,m}^+})} \calp
(M_{n,m})$:

\vskip 5pt

\nt {\bf Ree's Theorem} (see Theorem 2.1 of \cite{rees}) {\it If $I$
is generated by a regular sequence $f_1, \ldots, f_n$, then the map
$\phi: (S/I) [x_1,\ldots,x_n] \leftarrow \text{Gr}_I S$, sending
$x_i$ to the class $f_i$ in $I/I^2$ is an isomorphism. }

\vskip 5pt

Ree's Theorem thus implies that as vector spaces, we have $S = S/I
\otimes \Bbb C [f_1,\ldots, f_n]$, and in our context, $\calp
(M_{n,m}) \cong \calh_{n,m} \otimes \mathcal J_{n,m}$.

Thus we want to show that indeed $\{ r_{ij} \}$ form a regular
sequence in $I ( {\mathcal J_{n,m}^+})$. To show this, consider
the map from $M_{n,m}$ to the $m \times m$ symmetric matrices
$\calS^2 (\C^m)$ by putting the $r_{ij}$ in a matrix:
$$
Q: M_{n,m} \longrightarrow \calS^2 (\C^m) \qquad \text{where}\qquad Q(T)=T^tT, \qquad T\in M_{n,m}.
$$
First observe that this map is $O_n \times GL_m$ equivariant:
$$
Q(gTh)= h^t Q(T) h, \qquad g\in O_n, h \in GL_m.
$$
Further, the $O_n$ and $GL_m$ actions commute.

Let us study the fibers of the map $Q$.
Define the following rank $m$ matrices in $M_{n,m}$: For $k=0,1,\ldots, m$, let
$$
T_k =[ c_1 c_2 \ldots c_k c_{k+1} \ldots c_{m}] \in M_{n,m}
$$
where $\{ c_1, c_2, \ldots, c_k \}$ is an {\it orthonormal} set of
{\it non-isotropic} vectors in $\C^n$ and $\{ c_{k+1}, \ldots, c_m
\}$ is an {\it orthogonal} set of {\it isotropic} vectors in $\C^n$,
such that each of the $c_j$ with $j > k$ are orthogonal to the $c_j$
with $j \leq k$.  It is easy to see that
$$
Q(T_k) = \begin{bmatrix} I_k &0 \\ 0 &0 \end{bmatrix}
$$
where $I_k$ is the $k\times k$ identity matrix. Next, for
$i=0,1,\ldots,k$, define the sets\linebreak $\Theta_k \subset
M_{n,m}$ as follows:
$$
\Theta_k = \{ X \in M_{n,m} \mid Q(X) = Q(T_k) \text{ and } \rank \, X = \rank \, T_k =m \}.
$$
Let us remind our readers on the following version of Witt's Theorem (see Theorem 3.7.1 of \cite{howe-schur}):

\vskip 5pt

\nt {\bf Witt's Theorem} {\it Given two $n \times m$ matrices $T_1$ and $T_2$, there is an orthogonal $n \times n$ matrix $g$ such that $gT_1 = T_2$ if and only if $Q(T_1) = Q(T_2)$ and $\ker T_1 = \ker T_2$.
}

\vskip 5pt

Since $X$ and $T_k$ are of full rank, i.e., $\rank \, X = \rank \, T_k =m$ , thus $\ker X = \ker T_k = \{ 0 \}$. By Witt's Theorem, $\Theta_k$ is an $O_n$ orbit in the fiber $Q^{-1}(Q(T_k)) \subset M_{n,m}$.  The full rank condition gives the openness and denseness of this orbit.  The {\it null cone} (or null fiber) $NCQ = \Theta_0$ corresponds to $k=0$.

We claim that the fibers of the map $Q$ are all varieties of the same dimension, i.e., $Q$ is an {\it equi-dimensional} map:

\vskip 5pt

\nt {\bf Proposition} {\it Consider the $O_n \times GL_m$-equivariant map $Q$. Then for each $Y \in \calS^2 (\C^m)$,
\begin{enumerate}
\item[(a)] the fiber $Q^{-1}(Y)$ is an irreducible variety, invariant under $O_n$, and contains an open dense orbit of the form $\Theta_k \cdot h$, for some $h\in GL_m$ and some $k=0,1,\ldots,m$.
\item[(b)] the dimension of each fiber $Q^{-1}(Y)$ is $f = nm -\frac{m(m+1)}{2}$, which is independent of the fiber.
\end{enumerate}
}

\vskip 5pt

{\bf Remark:} The inverse image of any symmetric matrix contains
exactly one $O_n$ orbit consisting of invertible matrices. This is
dense in the whole inverse image. It is the condition that $n \geq
2m$ which permits injective maps into isotropic subspaces (i.e.,
such that the pulled-back form vanishes).

\nt {\bf Proof of Proposition:}  For $Y \in \calS^2 (\C^m)$, we can
find $h \in GL_m$ such that $h^tYh = \begin{bmatrix} I_k &0 \\ 0 &0
\end{bmatrix}$, for some $k=0,1,\ldots,m$.

If $X \in M_{n,m}$ is such that $Q(X)=Y$, then $Q(Xh) = h^tQ(X)h
=h^tYh = \begin{bmatrix} I_k &0 \\ 0 &0 \end{bmatrix}$. Further,
$\rank \, Xh = \rank \, X$. In other words, the fiber $Q^{-1}(Y)$
associated to $Y \in \calS^2 (\C^m)$ is the closure of an open dense
orbit given by $\Theta_k \cdot h$, and hence an irreducible variety.

The fact that the orbits $\Theta_k$ (and their translates) have the
same dimension follows from the computation of the dimension of the
pointwise stabilizer in $O_n$ of $T_k \in M_{n,m}$.  The pointwise
stabilizer can be easily computed for each $k =0,1,\ldots,m$. The
dimension of each fiber is $f = nm - \frac {m(m+1)}{2}$. (You can
see this dimension from the null cone pretty easily because the ring
of $O_n$ invariants is freely generated by $\frac{m(m+1)}{2}$
polynomials $r_{ij}$'s.) $\quad \square$

Since the mapping $Q$ has equi-dimensional fibres, if $V \subset
\calS^2 (\C^m)$ is any irreducible variety of dimension $e$, then
$Q^{-1} (V)$ will be an irreducible variety of dimension $e+f$. The
$Q(r_{ij})$'s are coordinates on $\calS^2 (\C^m)$, so the variety
defined by $d$ of them is a subspace of codimension $d$. It follows
that the pullback of this subspace by $Q$ is also irreducible and of
codimension $d$.  Therefore, the dimension of the zero set of
$r_{ij}$'s decrease by 1 at each stage, making $\{ r_{ij} \}$ a
regular sequence (see Lemma 4 on page 105 of \cite{matsumura}).
$\quad \square$

\vskip 10pt

\def\cprime{$'$} \def\cprime{$'$}

\end{document}